\documentclass{article}

\usepackage{xypic, amssymb, amsmath, amsthm}

\newtheorem{thm}{Theorem}[section]
\newtheorem{cor}[thm]{Corollary}
\newtheorem{lem}[thm]{Lemma}
\newtheorem{prop}[thm]{Proposition}

\newtheorem{claim}[thm]{Claim}
\theoremstyle{definition}
\newtheorem{defn}[thm]{Definition}
\theoremstyle{remark}
\newtheorem{rem}[thm]{Remark}
\numberwithin{equation}{section}

\newcommand{\set}[1]{\left\{#1\right\}}

\newcommand{\from}{\leftarrow}

\newcommand{\into}{\hookrightarrow}
\newcommand{\onto}{\twoheadrightarrow}
\newcommand{\id}{\mathbf{1}}
\newcommand{\quotes}[1]{\textquoteleft#1'}
\newcommand{\Ext}{\mbox{Ext}}
\newcommand{\RHom}{\mbox{RHom}}

\newcommand{\Hom}{\mbox{Hom}}
\newcommand{\End}{\mbox{End}}

\newcommand{\iso}{\stackrel{\sim}{\longrightarrow}}
\newcommand{\sslash}[1]{\,/\!/_#1\,}
\newcommand{\Tr}{\mbox{Tr}}

\begin{document}

\title{The $A_\infty$ Deformation Theory of a Point and the Derived Categories
of Local Calabi-Yaus}%
\author{Ed Segal\\
\small{Department of Mathematics}\\
\small{Imperial College London, SW7 2AZ, UK}\\
 \small{\textit{edward.segal@imperial.ac.uk}}
}%

\maketitle

\begin{abstract}

Let $A$ be an augmented algebra over a semi-simple algebra $S$. We
show that the Ext algebra of $S$ as an $A$-module, enriched with its
natural A-infinity structure, can be used to reconstruct the
completion of $A$ at the augmentation ideal. We use this technical
result to justify a calculation in the physics literature describing
algebras that are derived equivalent to certain non-compact
Calabi-Yau three-folds. Since the calculation produces
superpotentials for these algebras we also include some discussion
of superpotential algebras and their invariants.

\end{abstract}

\tableofcontents

\section{Introduction}\label{sectintro}

There are now various examples known of the phenomenon whereby a
variety $X$ can be derived equivalent to a non-commutative algebra
$A$. The pioneering example is due to Beilinson \cite{beilinson} who
proved that the derived category of $\mathbb{P}^n$ is generated by
the line bundles $\mathcal{O},...,\mathcal{O}(n)$. This equivalent to saying that
the functor
$$\RHom(\bigoplus_{i=0}^n \mathcal{O}(i), -): D^b(\mathbb{P}^n)
\to D^b(A) $$ is a derived equivalence between $\mathbb{P}^n$ and
the non-commutative algebra
$$A:=\End(\bigoplus_{i=0}^n \mathcal{O}(i))$$
In fact compact examples like this are rare, much more progress has been
made for non-compact examples, in particular for local models of resolutions
of singularities \cite{BKR},\cite{vandenbergh}.

The phenomenon is also well known in the physics literature. There
the variety $X$ should be a Calabi-Yau threefold, and we study the
type II superstring compactification on $X$ . Type B D-branes in the
theory correspond to objects in $D^b(X)$. It has been known since
the work of Douglas and Moore \cite{DM} that if a D-brane sits at
the centre of a singularity the effective theory on its world-volume
is a gauge theory whose content can be described by a quiver
diagram. This is the same as the mathematical results - the quiver
diagram is a presentation of an algebra $A$ which is derived
equivalent to a resolution of the singularity. Since then other
physical approaches (e.g \cite{HW}) have been found that produce an
effective quiver gauge theory from branes on $X$.

The example that we are interested is when $X = \omega$ is the
canonical bundle of a del Pezzo surface $Z$. This a Calabi-Yau
three-fold, and it is again \quotes{local} in that we may think of
it as the normal bundle to an embedded surface in a compact
Calabi-Yau.  A first step in describing $D^b(\omega)$ is to describe
$D^b(Z)$, and we specified that $Z$ should be a del Pezzo because in
that case Beilinson's approach has been generalised. What we do is
find a special collection of line bundles $\set{T_i}$ on $Z$ that
generate $D^b(Z)$, then as before $Z$ is derived equivalent to
$$A:= \End_Z(\oplus_i T_i)$$
As Bridgeland observed in \cite{bridgeland}, we have a similar
description of the derived category of $\omega$. If we pull up the
$T_i$ via the projection $\pi:\omega\to Z$ we find that they still
generate the derived category, so $\omega$ is derived equivalent
(under one further assumption) to
$$\tilde{A}:= \End_\omega(\oplus_i \pi^*T_i)$$
Both $A$ and $\tilde{A}$ can be presented as quiver algebras (with relations),
where the nodes
of the quiver correspond to the line bundles in the collection. Suppose we
have such a presentation of $A$. What do we have to do to it to produce a
presentation of $\tilde{A}$?

This was the question addressed, in rather more physical language,
by Aspinwall and Fidkowski in \cite{AF}. This paper is a
mathematical interpretation of their work, and of related physics
papers (\cite{AK}, \cite{lazaroiu}, \cite{BP} etc.). For the remainder of this
introduction we will discuss the answer to this question, leaving
out many subtleties and technicalities.

Suppose we have a presentation of the algebra $A$ as the path
algebra of a quiver $Q$ (with relations), where nodes of $Q$
correspond to the line bundles $T_i$. Then an $A$-module is
precisely a representation of the quiver that obeys the relations.
We have some obvious one-dimensional modules $S_i$ which are the
representations with just a one-dimensional vector space at the
$i$th node. The direct sum
$$\mathcal{S}=\bigoplus_i S_i$$ of these is a
representation which is one-dimensional at each node and with all
the arrows sent to zero maps.

If we pick projective resolutions of each $S_i$ then we can form the
dga
$$\RHom_A(\mathcal{S}, \mathcal{S})$$
and then, using the process of homological perturbation (\cite{GS}
etc...) transfer the dga structure to an $A_\infty$-structure on its
homology $\Ext_A(\mathcal{S},\mathcal{S})$.
Of course since $A$ is derived equivalent to $Z$ we could also view
the $S_i$ as being objects in $D^b(Z)$ and compute this
$A_\infty$-algebra there.

Now we consider the algebra $\tilde{A}$ corresponding to $\omega$.
This is also a quiver algebra on the same number of nodes, so has a
similar set of one-dimensional modules $\tilde{S_i}$. It is easy to
show that under the derived equivalence these map to the objects
$$\iota_*S_i \in D^b(\omega)$$
where $\iota: Z \to \omega$ is the zero section. We can again form
the sum
$$\tilde{\mathcal{S}}=\bigoplus_i \tilde{S}_i$$
 and the $A_\infty$-algebra
$$\Ext_{\tilde{A}}(\tilde{\mathcal{S}}, \tilde{\mathcal{S}})
= \Ext_{\omega}(\iota_*\mathcal{S},\iota_*\mathcal{S})$$

This new $A_\infty$ algebra has a straightforward relationship with
the previous one. By resolving the structure sheaf of the zero
section and using Serre duality on $Z$ one easily shows
\begin{equation}\label{cyccomp}\Ext_\omega(\iota_* \mathcal{S}, \iota_* \mathcal{S}) =
\Ext_Z(\mathcal{S},\mathcal{S}) \oplus
 \Ext_Z(\mathcal{S},\mathcal{S})[3]^\vee\end{equation}
 The two summands are dual under the
Calabi-Yau pairing on $D^b(\omega)$, and the $A_\infty$ structure
should be cyclic with respect to this pairing. In fact with a little
more work one can show that the $A_\infty$ structure is given by
formally extending the $A_\infty$ structure on
$\Ext_Z(\mathcal{S},\mathcal{S})$ to make it cyclic. We call this
procedure \textit{cyclic completion}.

Now we come to the key point:
\begin{claim} \label{claim} The algebra $A$ is determined by the
$A_\infty$-algebra $\Ext_A(\mathcal{S}, \mathcal{S})$, in that if
$\set{m_i}$ are the $A_\infty$-products on $\Ext_A(\mathcal{S},
\mathcal{S})$ then the map
$$(\oplus_i m_i)^\vee : \Ext^2_A(\mathcal{S}, \mathcal{S})^\vee \to
T^\bullet\Ext^1_A(\mathcal{S}, \mathcal{S})^\vee$$ is a presentation
of $A$. Similarly the $A_\infty$-algebra $$\Ext_{\tilde{A}}(
\tilde{\mathcal{S}}, \tilde{\mathcal{S}})$$ gives rise to a
presentation for $\tilde{A}$.
\end{claim}
This says that generators for $A$ are given by (the dual space to)
$\Ext^1_A( \mathcal{S}, \mathcal{S})$ and relations are given by
$\Ext^2_A( \mathcal{S}, \mathcal{S})$, with the form of the
relations being determined by the $A_\infty$ structure. If we split
$\mathcal{S}$ into its summands we see that this presentation is
actually of a quiver algebra: the generating arrows between nodes
$i$ and $j$ are given by $\Ext^1_A(S_i, S_j)$, and the relations on
paths between $i$ and $j$ are given by $\Ext^2_A(S_i, S_j)$.

This claim is the hard part of the argument, and Section \ref{secttheAinftydeftheoryofapt}
of this paper is devoted to the discussion and proof of it. However for now
we put it to one side and return to the question of determining
$\tilde{A}$.

Suppose that we have a presentation for $A$ of the form given in
Claim \ref{claim}. What is the corresponding presentation of
$\tilde{A}$? Using (\ref{cyccomp}):
$$\Ext^1_{\tilde{A}}(\tilde{S}_i, \tilde{S}_j) =
\Ext^1_A(S_i,S_j) \oplus
 \Ext^2_A(S_j,S_i)^\vee$$
 and
$$\Ext^2_{\tilde{A}}(\tilde{S_i}, \tilde{S_j}) =
\Ext^2_A(S_i,S_j) \oplus
 \Ext^1_A(S_j,S_i)^\vee$$
So the answer is that for each existing relation on paths from node
$j$ to node $i$ we should insert a new generating arrow going from
$i$ to $j$. Then for each existing generator going from $i$ to $j$
we put on one extra relation on paths going from $j$ to $i$. To
understand what the form of the relations should be we need to unpack our
definition of \quotes{cyclic completion}. It is easier to express the result
if we introduce the
notion of a \textit{superpotential}.

In fact from the physics perspective, working out the superpotential
is the primary goal, as it specifies the quiver gauge theory coming
from $\omega$. For the moment however we shall treat it just as the
following little trick from linear algebra. The spaces
$\Ext^1_{\tilde{A}}(\tilde{\mathcal{S}},\tilde{\mathcal{S}})$ and
$\Ext^2_{\tilde{A}}(\tilde{\mathcal{S}},\tilde{\mathcal{S}})$ are
dual under the Calabi-Yau pairing. Therefore the presentation
$$\Ext^2_{\tilde{A}}(\tilde{\mathcal{S}},\tilde{\mathcal{S}})^\vee \to
T^\bullet
\Ext^1_{\tilde{A}}(\tilde{\mathcal{S}},\tilde{\mathcal{S}})^\vee $$
alluded to in Claim \ref{claim} is given by an element
$$W \in T^\bullet
\Ext^1_{\tilde{A}}(\tilde{\mathcal{S}},\tilde{\mathcal{S}})^\vee$$
This is the superpotential for $\tilde{A}$. It is a formal
non-commutative polynomial in the generators, and taking partial
derivatives of it one recovers the relations. It is moreover
cyclicly symmetric since the $A_\infty$ structure is cyclic.

Now we can state the result. Suppose $A$ is given
by generators $\set{x_1,...,x_i}$ and relations
$\set{\rho_1,...,\rho_j}$ (which are formal expressions in the
$x_i$). Then the algebra $\tilde{A}$ is generated by the set
$\set{x_1,..,x_i,y_1,..., y_j}$ with relations coming from the
superpotential
$$W = \sum_{ \substack{ \mbox{\scriptsize cyclic} \\ \mbox{\scriptsize permutations}}}    \sum_j y_j\otimes \rho_j$$

\subsection{An example}

We illustrate the procedure with the prototypical example of
$\mathbb{P}^2$, with the line bundles $T_i = \mathcal{O}(i)$,
$i=0,1,2$. The endomorphism algebra $A$ of this collection is given
by the Beilinson quiver
$$ 
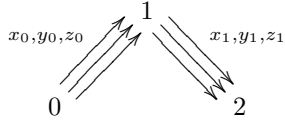
\xymatrix{
    & 1  \ar@<1ex>[dr]^{x_1,y_1,z_1}\ar[dr]\ar@<-1ex>[dr] &  \\
  0 \ar@<1ex>[ur]^{x_0,y_0,z_0}\ar[ur]\ar@<-1ex>[ur]&  & 2
   }$$
subject to the relations
$$\begin{aligned}
x_0y_1-y_0x_1&=0\\
y_0z_1-z_0y_1&=0\\
z_0x_1-x_0z_1&=0
\end{aligned}$$

Now we pass to the local Calabi-Yau $\omega =
\mathcal{O}(\mathbb{P}^2, -3)$, and pull up the line bundles. This
corresponds to cyclically completing the quiver algebra. Firstly we
insert extra arrows, dual to the relations. We have three relations,
each of which applies to paths from $T_0$ to $T_2$. Hence we should
insert three dual arrows from $T_2$ to $T_0$, so $\tilde{A}$ is
generated by the quiver
$$
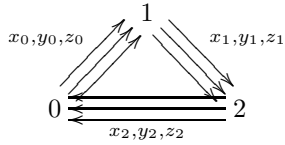
\xymatrix{
    & 1\ar@<1ex>[dr]^{x_1,y_1,z_1}\ar[dr]\ar@<-1ex>[dr] &  \\
  0 \ar@<1ex>[ur]^{x_0,y_0,z_0}\ar[ur]\ar@<-1ex>[ur]&
  & 2 \ar@<1ex>[ll]^{x_2,y_2,z_2}\ar[ll]\ar@<-1ex>[ll]
   }
$$
The superpotential is given by multiplying these new arrows by their
corresponding relations, so it is
\begin{eqnarray*}W &=& \sum_{ \substack{ \mbox{\tiny cyclic} \\ \mbox{\tiny permutations}}}
(x_0y_1- y_0x_1)z_2+ (y_0z_1- z_0y_1)x_2 + (z_0x_1-x_0z_1)y_2
\\
&=& \sum \epsilon^{ijk}x_i y_j z_k
\end{eqnarray*}
Now we compute the relations in $\tilde{A}$, which are given by
taking formal partial derivatives of $W$. Taking derivatives with
respect to the new generators just gives back the original three
relations. Taking derivatives with respect to the original
generators gives six new relations, each of which is a commutativity
relation of the form of the one of the original relations but lying
between a different pair of nodes.

According to our prescription the resulting algebra $\tilde{A}$
should be
$$\End_{\omega}(\pi^*\mathcal{O} \oplus \pi^*\mathcal{O}(1)\oplus\pi^*\mathcal{O}(2))$$
This is easily seen to be correct, since the latter is given by
$$
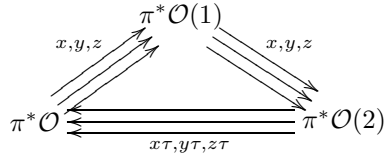
\xymatrix{
    & \pi^*\mathcal{O}(1)\ar@<1ex>[dr]^{x,y,z}\ar[dr]\ar@<-1ex>[dr] &  \\
  \pi^*\mathcal{O} \ar@<1ex>[ur]^{x,y,z}\ar[ur]\ar@<-1ex>[ur]&
  & \pi^*\mathcal{O}(2) \ar@<1ex>[ll]^{x\tau,y\tau,z\tau}\ar[ll]\ar@<-1ex>[ll]
   }
$$
where $\tau$ is the tautological section of $\pi^*\mathcal{O}(-3)$.

\subsection{The physical argument}\label{sectphysicalargument}

It is instructive to look at the physical arguments involved in
justifying Claim \ref{claim}. The set-up is type II superstring
theory on the ten-dimensional space $\omega\times \mathbb{R}^{3,1}$.
We have a D3-brane, which is a $(3+1)$-dimensional object, extending
in the flat directions, so from the point of view of $\omega$ it is
just a point $p$. The effective (i.e. low-energy limit) theory on
the world-volume of this brane is a gauge theory on
$\mathbb{R}^{3,1}$. The quiver diagram for $\tilde{A}$ specifies
this gauge theory - the nodes are $U(1)$ gauge groups, the arrows
are fields, and the relations are constraints on the fields.

In terms of the derived category this D3-brane is the skyscraper
sheaf $\mathcal{O}_p$. Under the derived equivalence between $Z$ and
$\tilde{A}$ this gets mapped to the $\tilde{A}$-module
$$\RHom_\omega(\bigoplus_i \pi^*T_i, \mathcal{O}_p) = \bigoplus_i (\pi^*T^\vee_i)|_p$$
This is a quiver representation that is one-dimensional at each
node.

The moduli space of $p$ is obviously just $\omega$. However, on the
other side of the derived equivalence it is also a moduli space
$\mathcal{M}$  of quiver representations that are one-dimensional at
each node, physically this is the vacuum moduli space of the quiver
gauge theory. We can construct this space as follows. Suppose we
have a presentation of $\tilde{A}$ as a quiver algebra with
generating arrows $V$ and some relations. Then a
$(1,...,1)$-dimensional representation is just an assignment of a
complex number to each generating arrow, such that the relations
hold. This means that the space of such representations is a
subvariety of $V^\vee$ cut out by the relations. Finally we must
quotient this space by the gauge action of
$\mathbb{C}^*\times...\times \mathbb{C}^*$ given by changing the
bases of the vector spaces at each node.

We now pick a K\"ahler metric on $\omega$, which gives us a notion
of stability for branes, and then deform the K\"ahler class to the
limit where the metric collapses the zero section $Z$ to a point. If our D3-brane
was sitting at a point $p\in Z$ then it becomes unstable in this
limit, and decays into a collection of so-called \textit{fractional}
branes, the $\tilde{S}_i$. We can see this mathematically in the
construction of $\mathcal{M}$. When we take the gauge group quotient
we should really pick a character $\chi$ of the gauge group and form
the GIT quotient $\mathcal{M}^\chi$. For appropriate characters this
should make the stable representations correspond precisely to
points $p\in \omega$, and thus $\mathcal{M}^\chi = \omega$. But if
we set $\chi=0$ then all representations corresponding to points in
$Z$ become semi-stable and S-equivalent to the origin in $V^\vee$,
which is the representation $\bigoplus_i \tilde{S}_i$. The moduli
space $\mathcal{M}^\chi$ is then the singularity obtained by
collapsing the zero section in $\omega$.

Now the physics of the D3-brane is encoded in the superpotential $W$
for the quiver gauge theory. This means that the equations of motion
for $p$ are the partial derivatives $\partial W$. However from the
construction of $\mathcal{M}$ we know that the equations restricting $p$
are precisely the relations in $\tilde{A}$, so in fact $W$ is a superpotential
in the mathematical sense for the algebra $\tilde{A}$.

On the other hand we can also see the behaviour of
$p$ by deforming $\oplus \tilde{S}_i$, since the deformation space
is just $\mathcal{M}$. These deformations will be governed by the
$A_\infty$-algebra
$$\Ext_{\tilde{A}}(\oplus \tilde{S}_i, \oplus \tilde{S}_i)$$
in the sense that if $W' \in T^\bullet (\Ext^1)^\vee$ encodes the
$A_\infty$ structure then the critical locus of $W'$ is the
deformation space of $\oplus_i \tilde{S}_i$. Thus $W'=W$ is the superpotential
for the quiver gauge theory, and hence for the algebra $\tilde{A}$.

\subsection{Notation and basics}

We will work over the ground field $\mathbb{C}$, although Section \ref{secttheAinftydeftheoryofapt}
works over an arbitrary ground field, and Section \ref{sectsuperpotentials} works over any
field of characteristic zero. $\mathbf{Alg}_\mathbb{C}$ is the category of
associative unital $\mathbb{C}$-algebras. Undecorated tensor products will be
over $\mathbb{C}$.

We will also need the category $\mathbb{C}^r$-$\mathbf{bimod}$ of
bimodules over the semi-simple ring $\mathbb{C}^r$. We denote the
obvious idempotents in $\mathbb{C}^r$ by $1_1,...,1_r$, then any
\linebreak$V\in\mathbb{C}^r$-$\mathbf{bimod}$ is a direct sum of the
subspaces
$$V_{ij} := 1_i.V.1_j$$
We may think of $V$ as a \quotes{categorified} vector space.

Let $\mathbf{Alg}_\mathbb{C}^r$ be the category of $\mathbb{C}^r$-algebras, i.e.
associative unital algebra objects in $\mathbb{C}^r$-$\mathbf{bimod}$.
Equivalently this is the category of $\mathbb{C}$-linear categories whose
objects form an ordered set of size $r$, if we only allow functors
that preserve the ordering on the objects.

Any algebra $A\in \mathbf{Alg}_\mathbb{C}^r$ may be pictured as a quiver
algebra (with relations) - just pick a basis for each $A_{ij}$, then
$A$ is a quotient of the path-algebra of the quiver with $r$ nodes
and arrows given by the basis elements. We may also consider $A$ as
an object of $\mathbf{Alg}_\mathbb{C}$ equipped with an ordered complete set
of orthogonal idempotents $\set{1_1,...,1_r}$.

If $V$ is any $\mathbb{C}^r$-bimodule then it generates a free $\mathbb{C}^r$-algebra
$$ TV : = \bigoplus_n V^{\otimes_{_{\mathbb{C}^r}} \, n}$$
and a completed algebra
$$ \hat{T} V : = \prod_n  V^{\otimes_{_{\mathbb{C}^r}} \, n}$$

$\mathbf{Alg}_\mathbb{C}^r$ admits a symmetric monoidal product
$\underline{\otimes}$
 given by
$$(A\underline{\otimes} B)_{ij} = A_{ij}\otimes B_{ij} $$

Note that this is certainly not $A\otimes_{\mathbb{C}^r} B$, in general
$A\otimes_{\mathbb{C}^r} B$ does not have an algebra structure.

An \textit{augmentation} of an algebra $A\in \mathbf{Alg}_\mathbb{C}^r$ is a
splitting $p: A \to \mathbb{C}^r$ of the inclusion $\mathbb{C}^r \into A$ of the
identity arrows, or equivalently a choice of a two-sided ideal
$\bar{A}\subset A$ such that $A/\bar{A} = \mathbb{C}^r$. Alternatively we may
think of $A$ as an algebra in $\mathbf{Alg}_\mathbb{C}$ for which we have
chosen $r$ $\mathbb{C}$-points $p: A \to \mathbb{C}^r$ and then split $p$. We denote
the category of augmented algebras by $\mathbf{Alg}_\mathbb{C}^{\star r}$.
Morphisms must respect the augmentations.

A module always means a left module. If we are picturing $A\in
\mathbf{Alg}_\mathbb{C}^r$ as a quiver algebra then a module over $A$ is
precisely a representation of the quiver (that respects the
relations). It is also the same as a functor $A \to \mathbf{Vect}$.

\subsection{Acknowledgements}

Thanks to Kai Behrend, Lieven Le Bruyn, Tom Coates, Alessio Corti,
Joel Fine, Dominic Joyce, Alistair King, Rapha\"el Rouquier, Jim
Stasheff and Bal\'azs Szendr\"oi for many helpful conversations,
comments and ideas. Particular thanks are due to Kevin Costello and
Paul Seidel for their help and hospitality at the University of
Chicago, and to Tom Bridgeland and Simon Donaldson who examined this material
when it was presented as my doctoral thesis. Finally I owe an enormous debt to my supervisor Richard
Thomas, for teaching me (amongst other things) that having
 a concrete
example sometimes means more than just specifying the ground field
to be $\mathbb{C}$.

\section{The $A_\infty$ Deformation Theory of a Point}
\label{secttheAinftydeftheoryofapt}

In this section we address the following claim, which we made in the introduction:
suppose we have an appropriate set $\set{S_i}$ of one-dimensional modules for some algebra $A$.
Then we can reconstruct $A$ from the $A_\infty$-algebra $\Ext_A(\oplus
S_i, \oplus S_i)$.

In fact if we assume that $A$ is graded, and that $A_0=\oplus S_i$, then this statement has been part of the mathematical folklore for some
time. The result is claimed (although not
proven) by Keller \cite{keller2} for a particular class of graded
algebras, and his statement is closely related to a result of Laudal
\cite{laudal}, who uses the terminology of Massey products. The
fullest investigation to date appears to be the work of Lu,
Palmieri, Wu and Zhang \cite{LPWZ}.

Let us start by explaining the statement a little. Let $A$ be an
$\mathbb{N}$-graded algebra over $\mathbb{C}$, and for simplicity let
$A_0 =\mathbb{C}$. Now suppose we are given a presentation
$$ A = TV / (\iota R)$$
so $A$ is generated by a vector space $V$, modulo the two-sided
ideal generated by a space of relations $R$ under an inclusion
$$\iota: R \to TV $$
Assume that the presentation is minimal, in the sense that $V$ and
$R$ are of minimal dimension. Then using the free resolution
$$ ... \to A\otimes R \to A\otimes V \to A \to \mathbb{C} \to 0 $$
of $\mathbb{C} = A_0$ it is elementary to show that $V$ must be dual to
$\Ext^1_A(\mathbb{C},\mathbb{C})$, and $R$ must be dual to $\Ext^2_A(\mathbb{C},\mathbb{C})$. Hence we
might ask: if we are just given $\Ext^\bullet_A(\mathbb{C},\mathbb{C})$, can we
recover $A$?

We know immediately that $A$ is generated by the space
$V:=(\Ext^1_A(\mathbb{C},\mathbb{C}))^\vee$, and that relations are counted by the
space $R:=(\Ext^2_A(\mathbb{C},\mathbb{C}))^\vee$, but we still need to know what form
these relations take, i.e. we need the map
$$\iota: R \to T V $$
or dually, a map
$$\iota^\vee: \hat{T}\Ext^1_A(\mathbb{C},\mathbb{C}) \to \Ext_A^2(\mathbb{C},\mathbb{C}) $$
We certainly have something that might be a part of this map, namely
the usual Yoneda (wedge) product, which is a map
$$ \Ext^1_A(\mathbb{C},\mathbb{C})^{\otimes 2} \to \Ext_A^2(\mathbb{C},\mathbb{C}) $$
If we knew that our relations were purely quadratic then we might
reasonably conjecture that this dualising this map gave a
presentation of $A$. In fact although this is true for many algebras
it is false in general - the study of those algebras for which it
works is the subject of classical Koszul duality. What happens when
our relations are definitely not just quadratic? Then we would need,
in addition to the bilinear Yoneda product, some \quotes{higher}
multi-linear products
$$m_i: \Ext^1_A(\mathbb{C},\mathbb{C})^{\otimes i} \to \Ext_A^2(\mathbb{C},\mathbb{C})$$
Fortunately these higher products do exist (though not quite
canonically), they form an $A_\infty$-structure on
$\Ext^\bullet_A(\mathbb{C},\mathbb{C})$ which measures the failure of the dga
$\RHom_A(\mathbb{C},\mathbb{C})$ to be formal. Furthermore when you dualize they do
indeed give a presentation of $A$. It is this result (essentially
our Theorem \ref{gradedresult}) that is proven by \cite{LPWZ}.

One of our original aims was to prove this result for the case $A_0
= \mathbb{C}^r$, i.e. when $A$ is a graded quiver algebra (with relations) on
$r$ vertices. However, given the proof in \cite{LPWZ} this is easy -
you simply change your ground category from vector spaces to the
category of $\mathbb{C}^r$-bimodules (which one may picture as vector spaces
strung between $r$ vertices) and the same proof works. Instead we
take a different tack which we feel is a bit more conceptual.

It seemed to us that the graded hypothesis was a little unnatural.
We instead ask what happens if we take an arbitrary algebra $A$ with
a one-dimensional module $S$ and perform the same construction, i.e.
take the Yoneda algebra $\Ext_A(S,S)$ equipped with $A_\infty$
products $\set{m_i}$, then dualize the map
$$m = \oplus_i m_i : T \Ext^1_A(S,S) \to \Ext^2_A(S,S)$$
to get the presentation of a new algebra
$$ E := \frac{\hat{T} \Ext^1_A(S,S)^\vee}{( m^\vee \Ext^2_A(S,S)^\vee)} $$
What is this new algebra? Firstly note that a one-dimensional module
is just a map $p: A \to \mathbb{C}$. Hence if $A$ is commutative then this is
simply a closed point of the affine scheme Spec$(A)$, and the module
is its sky-scraper sheaf. It is then geometrically obvious that the
algebra $E$ can only depend on a formal neighbourhood of the point
$p$. In fact the result is that $E$ is precisely the formal
neighbourhood of $p$, i.e. it is the completion of $A$ at the kernel
of $p$. We explain this result (which contains nothing new)
informally in Section \ref{sectsketch}, the key point is that the
$A_\infty$-algebra $\Ext_A(S,S)$ controls the deformations of the
module $S$ and hence those of $p$.

This is of course the case $r=1$, in general we wish to pick $r$
points
$$p =\oplus p_i:A \to \mathbb{C}^r$$
 so that $\mathbb{C}^r$ becomes an
$A$-module (strictly speaking we must also
choose a splitting of the the map $p$, so this is more like choosing a single point
of a $\mathbb{C}^r$-algebra). Then our main result
(Theorem \ref{hatE=hatA}) is that performing the above construction
on $\Ext_A(\mathbb{C}^r, \mathbb{C}^r)$ again produces the completion of $A$ at the
kernel of $p$.

If we stick with a commutative $A$ then this generalization is
trivial, since $\Ext_A(\mathbb{C}^r, \mathbb{C}^r)$ splits as a direct product over
the different points. If $A$ is non-commutative however this is no
longer true and the proof becomes rather more difficult. In
particular it is not correct to study deformations of $\mathbb{C}^r$ as an
$A$-module, instead one should follow \cite{laudal} and study
non-commutative deformations of the category whose objects are these
$r$ one-dimensional modules. The technical challenge of this paper
is checking that the sketch proof given for the $r=1$ commutative
case continues to work in the general setting, which means firstly
checking that the non-commutative deformations of a set of modules
are governed by the $A_\infty$-category of their $\Ext$ groups
(Section \ref{sectdefthyofsetsofmods}) and secondly relating
deformations of a set of one-dimensional modules to deformations of
the corresponding set of points (Section
\ref{sectdefingasetofpoints}).

Having understood this more general situation it is then
straightforward to deduce the required result for graded algebras,
which we do in Section \ref{sectthegradedcase}. This is because the
completion of an $\mathbb{N}$-graded algebra at its
positively-graded ideal contains the original algebra in a natural
way.

\subsection{A geometric sketch} \label{sectsketch}

Let $X =$ Spec $A$ be an affine scheme over $\mathbb{C}$, and let $p:A\to \mathbb{C}$
be a point. Obviously the formal deformations of $p$ see
precisely a formal neighbourhood of $p$ in $X$, the algebraic way to
say this is that the formal deformation functor of $p$ is
pro-represented by the completion $\hat{A}_p$ of $A$ at the kernel
of $p$.

It is easy to show that the deformation theory of $p$ is precisely
the same as the deformation theory of the associated
\quotes{sky-scraper} sheaf $\mathcal{O}_p$, i.e. the 1-dimensional
$A$-module given by $p$. In accordance with the philosophy of dga
(or dgla, see Remark \ref{rem}) deformation theory, the deformations
of $\mathcal{O}_p$ are governed by the differential graded algebra
$\RHom_X(\mathcal{O}_p,\mathcal{O}_p)$, by which we mean that formal
deformations of $\mathcal{O}_p$ are formal solutions of the
\textit{Maurer-Cartan} equation
$$MC : \RHom^1_X(\mathcal{O}_p,\mathcal{O}_p) \to \RHom^2_X(\mathcal{O}_p,\mathcal{O}_p)$$
$$ MC(a):= da + a^2 = 0$$
taken modulo the \quotes{infinitesimal gauge action} of
$\RHom^0_X(\mathcal{O}_p,\mathcal{O}_p)$. This resulting
\quotes{formal deformation space} is, as we just said, simply the
formal scheme $\hat{A}_p$.

According to Kontsevich \cite{kontsevich} the formal deformation
theory attached to a dga is a homotopy invariant, so we may replace
our dga by any quasi-isomorphic $A_\infty$-algebra and compute the
deformations there instead. The Maurer-Cartan equation picks up
higher terms from the $A_\infty$ structure and becomes the
\textit{Homotopy Maurer-Cartan equation}:
$$HMC(a):=\sum_i m_i(a^{\otimes i}) =0$$
and there are similar homotopy corrections to the gauge action. In
particular, using the process of homological perturbation, we may
replace $\RHom_X(\mathcal{O}_p,\mathcal{O}_p)$ by its homology
$\Ext_X(\mathcal{O}_p,\mathcal{O}_p)$ equipped with an appropriate
$A_\infty$-structure.

Since $\Ext^0_X(\mathcal{O}_p,\mathcal{O}_p)=\mathbb{C}$, the gauge action is
now trivial, so the formal deformation space is the formal zero
locus of
$$HMC: \Ext^1_X(\mathcal{O}_p,\mathcal{O}_p) \to \Ext^2_X(\mathcal{O}_p,\mathcal{O}_p)$$
The algebra of functions on this formal scheme is the formal power
series ring on $\Ext^1$ modulo the ideal generated by the $\mathbb{C}$-linear
dual of $HMC$, so we have shown
 $$\frac{\mathbb{C}[[\;\Ext^1_X(\mathcal{O}_p,\mathcal{O}_p))^\vee\;]]}
 { \left( \;HMC^\vee(\Ext^2_X(\mathcal{O}_p,\mathcal{O}_p)^\vee)\;\right)} = \hat{A}_p$$

This says that formally around $p$, $X$ is cut out of $T_p X$ by the
HMC equation, with $\Ext^2_X(\mathcal{O}_p,\mathcal{O}_p)$ being a
canonical space of obstructions. Of course if $p$ is a smooth point the denominator
of this expression should vanish, since there are no (commutative) obstructions.
However in that case $\Ext^2$ is precisely the commutativity relations, which
tells us that we are really measuring non-commutative obstructions, and that the numerator should really be non-commutative power series.

\begin{rem}\label{rem}It is more traditional to control deformations
with dg-Lie (or $L_\infty$) algebras, but in this paper we will always in
fact have a dg (or
$A_\infty$) algebra. We should really take the associated commutator
algebra, as this is all that the deformation theory depends on, but we shall
not bother to do so.\end{rem}

\subsection{Deformation theory of sets of modules}
\label{sectdefthyofsetsofmods}

Let $A\in \mathbf{Alg}_\mathbb{C}$ and let
$\mathcal{M}=\set{M_1,...,M_r}$ be a set of $A$-modules. We show how
the non-commutative deformation theory of $\mathcal{M}$ as developed
by Laudal \cite{laudal} may be viewed as a dga deformation problem.

If we wanted to deform a single module $M$ then we would just deform
the module map
$$ \mu: A \to \End_\mathbb{C}(M)$$
When we have a set of modules we can form the endomorphism algebra
$$\End_\mathbb{C}(\mathcal{M}) := \End_\mathbb{C}(\oplus_i M_i)$$
and we could deform the map
$$ \mu =\oplus_i \mu_i: A \to \End_\mathbb{C}(\mathcal{M})$$
If we just treat this as a map of $\mathbb{C}$-algebras and deform
it then we are just studying deformations of $\oplus_i M_i$ as an
$A$-module. We wish to do something slightly different, and use the
fact that $\End_\mathbb{C}(\mathcal{M})$ is actually a
$\mathbb{C}^r$-algebra.

Recall that $\mathbf{Alg}_\mathbb{C}^{\star r}$ is the category of
augmented algebras over $\mathbb{C}^r$. Let
$$\mathbf{Art}_\mathbb{C}^r \subset \mathbf{Alg}_\mathbb{C}^{\star
r}$$ be the subcategory of consisting of algebras $(R,\mathfrak{m})$
for which the augmentation ideal $\mathfrak{m}$ is nilpotent. Of
course $\mathbf{Art}_\mathbb{C}^1$ is just the category of Artinian
local (non-commutative) $\mathbb{C}$-algebras. Recall also the
product of two $\mathbb{C}^r$-algebras that we defined by
$$(A\underline{\otimes} B)_{ij} = A_{ij}\otimes B_{ij} $$

\begin{defn}\label{defndefofmoduleset}\cite{laudal} For $(R,\mathfrak{m})\in \mathbf{Art_\mathbb{C}^r}$, an $R$-deformation
of $\mathcal{M}$ is a map of $\mathbb{C}$-algebras
$$\mu_R: A \to \End_\mathbb{C}(\mathcal{M})\underline{\otimes} R $$
which reduces modulo $\mathfrak{m}$ to the given module maps
$$\bigoplus_i  \mu_i : A \to \bigoplus_i \Hom_\mathbb{C}(M_i,M_i) $$
Two $R$-deformations are equivalent if they differ by an inner
automorphism of $\mathbb{C}$-algebras in
$\End_\mathbb{C}(\mathcal{M})\underline{\otimes} R$.
\end{defn}

Let $\mathcal{D}ef_\mathcal{M}: \mathbf{Art}_\mathbb{C}^r \to
\mathbf{Set}$ be the resulting deformation functor.

We are going to present this as a dga (and later $A_\infty$)
deformation problem. It is well known that the deformation functor
of a single module $M$ can be seen as a dga deformation problem - it
is controlled by the dga $\RHom_A(M,M)$. We now show that a similar
statement is true for our deformation functor
$\mathcal{D}ef_\mathcal{M}$, but since we are deforming
$\mathbb{C}^r$-algebras we look not for a dga but for a dg-category
with $r$ objects.

\begin{defn} Let $(\mathcal{A}^\bullet, d, m)$ be a dga over $\mathbb{C}^r$. The \textit{deformation functor}
associated to $\mathcal{A}$ is the functor
$$\mathcal{D}ef_{\mathcal{A}}: \mathbf{Art}_\mathbb{C}^r \to \mathbf{Set}$$
which sends $(R, \mathfrak{m})$ to the set
$$\set{a\in \mathcal{A}^1\underline{\otimes} \mathfrak{m} ; \;\; da +
 m(a\otimes a)= 0} / \sim $$
where the equivalence relation $\sim$ is given by taking the
exponential of the following action of the commutator Lie algebra of
$\mathcal{A}^0\underline{\otimes}\mathfrak{m}$
$$ b: a \to a + db - [b,a]$$
\end{defn}

For our deformation problem the obvious choice of
dg-$\mathbb{C}^r$-algebra is
$$\RHom_A(\mathcal{M},\mathcal{M}) = \bigoplus_{i,j} \RHom_A(M_i, M_j)$$
This is only defined up to quasi-isomorphism. To produce models for
it we need to resolve each $M_i$, which we may do using the
following standard construction:

\begin{defn} For an $A$-module $M$ the \textit{bar
resolution} of $M$ to be the complex of free $A$-modules
(concentrated in non-positive degrees)
$$B(A,M)^{-t} := A^{\otimes t+1}\otimes M$$
with differential given by
$$ d(a_1\otimes...\otimes a_t\otimes m) = \sum_{s=2}^t
(-1)^s a_1\otimes...\otimes a_{s-1}a_s\otimes ... \otimes a_t\otimes
m \;\;-\;\; (-1)^t a_1\otimes...\otimes a_t m $$
\end{defn}

\begin{lem} The module map $\mu: A\otimes M \to
M$ induces a quasi-isomorphism
$$\mu: B(A,M) \to M$$
\end{lem}
\begin{proof}
Since $A$ is unital $\mu$ is a surjection, and $B(A,M)$ is acyclic
in all negative degrees since
$$ d(1_A\otimes \mathbf{b}) = \mathbf{b} - 1_A\otimes d(\mathbf{b})$$
for any $\mathbf{b} \in B(A,M)^{<0}$.\end{proof}

Hence one model for $\RHom_A(\mathcal{M},\mathcal{M})$ is given by
the dg-category $\mathcal{E}$ whose hom-sets are
$$ \mathcal{E}_{ij}:= \mbox{Hom}_A(B(A,M_i), B(A, M_j))$$
However there is a simpler candidate. Consider the dg-category
$\mathcal{H}$ whose hom-sets are
$$\mathcal{H}_{ij} := \mbox{Hom}_A(B(A,M_i), M_j) $$
with composition
$$ (f\bullet g)(a_1\otimes...\otimes a_{s+t+1}\otimes m) :=
f(a_1\otimes...\otimes a_{s+1}\otimes g( 1_A\otimes
a_{s+2}\otimes...\otimes a_{s+t+1}\otimes m)) $$ for homogeneous
maps $f,g$ of degrees $s$ and $t$. This composition was obtained as
follows: $B(A,\mathbb{C})$ is naturally a coalgebra under the
\quotes{shuffle} coproduct, and $B(A,\mathcal{M})$ is a comodule
over it. We are letting $f\bullet g = f(\id\otimes g)\mu$ where
$\mu$ is the comodule map. Now for each $i,j$ we have a
quasi-isomorphism
$$\mu_j: B(A,M_j) \to M_j$$
which induces a quasi-isomorphism of chain complexes
$$ \mu_j\circ: \mathcal{E}_{ij} \to \mathcal{H}_{ij} $$
(since the $B(A,M)^i$ are free). These do not form a map of
dg-categories since they do not respect the compositions, but they
do have a right-inverse which is a map of dg-categories:

\begin{lem}\label{H=E} Let $$\Psi: \mathcal{H}^t \to \mathcal{E}^t$$
be given by
$$ \Psi(f)(a_1\otimes...\otimes a_{s+1}\otimes m) :=
a_1\otimes...\otimes a_{s-t+1} \otimes f( 1_A \otimes
a_{s-t+2}\otimes..\otimes a_{s+1}\otimes m)  $$ for $s\geq t$ and
$$ \Psi(f)(a_1\otimes...\otimes a_{s+1}\otimes m) := 0$$
 for $s<t$. Then $\Psi$ is a map of dg-categories such that
$(\mu\circ)\Psi = \id_\mathcal{H}$, hence it is a quasi-isomorphism
of dg-categories.
\end{lem}
\begin{proof} Elementary (though tedious) from definitions. \end{proof}

We use this model $\mathcal{H}$ for
$\RHom_A(\mathcal{M},\mathcal{M})$ to show that this dg-category
controls the deformations of $\mathcal{M}$, at least (as in the
single module case) up to module automorphisms of the $M_i$.

\begin{prop} \label{defMcontrolledbyH} $\mathcal{D}ef_\mathcal{M}$ is a
 quotient by $\mbox{Aut}_A(\bigoplus_i M_i)$ of the deformation functor
 $\mathcal{D}ef_{\mathcal{H}}$ associated
to the dg-category $\mathcal{H}$.
\end{prop}
\begin{proof}
Let $R$ be an object of $\mathbf{Art}_\mathbb{C}^r$. Using the
splitting $R = \mathbb{C}^r \oplus \mathfrak{m}$ we can write any
$R$-deformation $\mu_R$ of $M$ as
$$\mu_R = \bigoplus_i \mu_i + \tilde{\mu}_R$$
where $\tilde{\mu}_R :A \to
\End_\mathbb{C}(\mathcal{M})\underline{\otimes}\mathfrak{m}$. The
$R$-points of $\mathcal{D}ef_\mathcal{H}$ are those (equivalence
classes of) elements of $\mathcal{H}^1\underline{\otimes}
\mathfrak{m}$ that obey the Maurer-Cartan equation. However,
\begin{eqnarray*} \mathcal{H}^1\underline{\otimes} \mathfrak{m} &=&
\bigoplus_{i,j} \Hom_A(A^{\otimes 2}\otimes M_i, M_j)\otimes
\mathfrak{m}_{ij}\\
&=& \Hom_\mathbb{C}(A, \bigoplus_{i,j}\Hom_\mathbb{C}(M_i,
M_j)\otimes
\mathfrak{m}_{ij})\\
&=& \Hom_\mathbb{C}(A, \End_\mathbb{C}(\mathcal{M})
\underline{\otimes} \mathfrak{m} )
\end{eqnarray*}
and for an element $\tilde{\mu}_R\in
\mathcal{H}^1\underline{\otimes}\mathfrak{m}$ the Maurer-Cartan
equation is precisely the condition that $\bigoplus_i \mu_i +
\tilde{\mu}_R$ is a map of algebras.

Now we compare the equivalence relations on each side. The class of
$\mu_R = \bigoplus_i \mu_i\label{dfg} + \tilde{\mu}_R$ in
$\mathcal{D}ef_\mathcal{M}(R)$ is its orbit under conjugation by the
the subgroup
$$\mbox{Stab}(\bigoplus_i \mu_i) $$ of the group of invertible elements in
$\End_\mathbb{C}(\mathcal{M})\underline{\otimes} R$. We have the
obvious factorization
$$  (1 + \End_\mathbb{C}(\mathcal{M})\underline{\otimes} \mathfrak{m}) \;\; \to\;\; \mbox{Stab}(\bigoplus_i \mu_i)
 \;\; \;\;\to\;\; \mbox{Aut}_A(\bigoplus_i M_i) $$
The Lie algebra of $(1 + \End_\mathbb{C}(\mathcal{M})$ is the
commutator algebra of
$$\End_\mathbb{C}(\mathcal{M})\underline{\otimes} \mathfrak{m} =
\mathcal{H}^0\underline{\otimes} \mathfrak{m}$$ The equivalence
relation on $\mathcal{D}ef_{\mathcal{H}}(R)$ is given by integrating
the \quotes{infinitesimal gauge action} of
$\mathcal{H}^0\underline{\otimes} \mathfrak{m}$, but this action is
precisely the derivative of conjugation. Hence the orbits in
$\mathcal{D}ef_\mathcal{M}(R)$ are the quotients of the orbits in
$\mathcal{D}ef_\mathcal{H}(R)$ under the residual action by
$\mbox{Aut}_A(\bigoplus_i M_i)$.
\end{proof}

\begin{cor}  If the $M_i$ are simple and distinct then
$\mathcal{D}ef_\mathcal{M} = \mathcal{D}ef_\mathcal{H}$
\end{cor}

In light of Lemma \ref{H=E} the homology of $\mathcal{H}$ is the
category $$\Ext_A(\mathcal{M})=\bigoplus_{i,j}\Ext_A(M_i,M_j)$$
Using homological perturbation, we may put an $A_\infty$-structure
on this category (unique up to $A_\infty$-isomorphism) such that it
is $A_\infty$-quasi-isomorphic to $\mathcal{H}$. We can now use this
$A_\infty$-algebra to compute the deformation functor of
$\mathcal{M}$. \pagebreak

\begin{defn} Let $(\mathcal{A}^\bullet, m_i)$ be an $A_\infty$-algebra over $\mathbb{C}^r$. The \textit{deformation functor}
associated to $\mathcal{A}$ is the functor
$$\mathcal{D}ef_{\mathcal{A}}: \mathbf{Art}_\mathbb{C}^r \to \mathbf{Set}$$
which sends $(R, \mathfrak{m})$ to the set
$$\set{a\in \mathcal{A}^1\underline{\otimes} \mathfrak{m} ; \;\; \sum_i m_i(a)= 0} / \sim $$
The equivalence relation $\sim$ is generated by the following map
from $\mathcal{A}^0\underline{\otimes}\mathfrak{m}$ to vector fields
on $\mathcal{A}^1\underline{\otimes} \mathfrak{m}$:
$$ b: a \to a + \sum_{n\geq 1} (-1)^{n(n+1)/2}\frac1n
\sum_{t=0}^{n-1}(-1)^t m_n(a^{\otimes t}\otimes b \otimes a^{\otimes
n-t-1}) $$
\end{defn}

We have obtained this from the usual definition of the deformation
functor of an $L_\infty$-algebra (e.g. \cite{lazaroiu}) in the case
that the $L_\infty$-algebra is actually the commutator algebra of an
$A_\infty$-algebra.

\begin{cor} \label{Def_McontrolledbyExt} $\mathcal{D}ef_\mathcal{M}$ is a quotient by
 $\mbox{Aut}_A(\bigoplus_i M_i)$ of the deformation functor
associated to the $A_\infty$-category $\Ext_A(\mathcal{M})$.
\end{cor}
\begin{proof}
It is standard (e.g. \cite{kontsevich}) that for $A_\infty$-algebras
over a field the usual commutative deformation functor is a homotopy
invariant. The key point of the standard proof is that any
$A_\infty$-algebra is isomorphic to the direct sum of a minimal one
and a linear contractible one, but in fact this holds when the base
category is any semi-simple linear monoidal category \cite{lefevre}
so it works over $\mathbb{C}^r$. The remainder of the proof consists
of checking three things: that $A_\infty$-morphisms induce natural
transformations of deformation functors, that deformation functors
commute with direct sums, and that the deformation functor
associated to a linear contractible $A_\infty$-algebra is trivial.
These are easily checked to hold for our non-commutative deformation
functors as well.
\end{proof}

\subsection{Deforming a set of points}\label{sectdefingasetofpoints}

Let $A$ be a $\mathbb{C}$-algebra, and let $p:A \to \mathbb{C}^r$ be
a set of $r$ $\mathbb{C}$-points of $A$. Each point $p_i: A \to
\mathbb{C}$ gives a one-dimensional $A$-module which we call
$\mathcal{O}_{p_i}$. We wish to relate the deformations of the
points $p$ to the deformations of the set of modules $\mathcal{M}:=
\set{\mathcal{O}_{p_i}}$.

In fact to achieve this in the way that we want we have to start
with a little more data - we have to choose a splitting of the map
$p$, so that $A$ becomes an augmented $\mathbb{C}^r$-algebra. Then
we can consider deformations of the map $p$ in the category
$\mathbf{Alg}_\mathbb{C}^r$, these are just $\mathbb{C}^r$-algebra
maps
$$p_R : A \to R $$
where $R\in\mathbf{Art}_\mathbb{C}^r$ and $p_R$ reduces to $p$
modulo the augmentation ideal in $R$. In other words $p_R$ must be a
morphism of augmented $\mathbb{C}^r$-algebras, so the deformation
functor of $p$ (in $\mathbf{Alg}_\mathbb{C}^r$) is
$$ \mathbf{Alg}_\mathbb{C}^{\star r}(A, -): \mathbf{Art}_\mathbb{C}^r \to \mathbf{Set} $$
How does this relate to the deformation functor
$\mathcal{D}ef_\mathcal{M}$ of the set of the modules?

\begin{prop} \label{defMrepbyA}
For any $(R,\mathfrak{m})\in\mathbf{Art}_\mathbb{C}^r$ there is a
functorial isomorphism
 $$\mathcal{D}ef_\mathcal{M}(R) = \mathbf{Alg}_\mathbb{C}^{\star r}(A,R) \;/\; \set{\mbox{Inner $\mathbb{C}^r$-algebra automorphisms}}$$
\end{prop}
\begin{proof}
Recall (Definition \ref{defndefofmoduleset}) that
$\mathcal{D}ef_\mathcal{M}(R)$ is the set of $\mathbb{C}$-algebra
maps
$$p_R: A \to  \End_\mathbb{C}(\mathcal{M})\underline{\otimes}R$$
that reduce to $p$ modulo $\mathfrak{m}$, taken up to inner
automorphism of $\mathbb{C}$-algebras. Now the hom-sets of
$\End_\mathbb{C}(\mathcal{M})$ are all 1-dimensional, so
$$ \End_\mathbb{C}(\mathcal{M})\underline{\otimes}R = R$$
for any $R$ (this fact is the reason that we defined our
deformations in terms of the product $\underline{\otimes}$). So to
prove the proposition we just need to check that
$\mathbb{C}$-algebra maps $A \to R$ that preserve the augmentations
are the same thing as augmented $\mathbb{C}^r$-algebra maps, once we
have quotiented out by inner automorphisms on both sides. This is
done in the following lemma.
\end{proof}

The difference between an inner $\mathbb{C}$-algebra automorphism of
$R$ and an inner $\mathbb{C}^r$-algebra automorphism is that the
former is conjugation by an arbitrary element $r\in R$ whereas for
the latter we must take $r$ from the diagonal subalgebra
$$\bigoplus_i R_{ii} \subset R $$
in order to preserve the $\mathbb{C}^r$-algebra structure.

\begin{lem} Let $(B,p)$ be an augmented $\mathbb{C}^r$-algebra in $\mathbf{Alg}_\mathbb{C}^{\star r}$ and
let\linebreak $(R,\mathfrak{m})\in\mathbf{Art}_\mathbb{C}^r$. Then
maps of $\mathbb{C}$-algebras $B \to R$ preserving the
augmentations, taken up to inner automorphism of
$\mathbb{C}$-algebras, biject with maps of augmented
$\mathbb{C}^r$-algebras $B\to R$ taken up to inner automorphism of
$\mathbb{C}^r$-algebras.
\end{lem}
\begin{proof} Obviously $\mathbb{C}^r$-algebra maps form a subset of $\mathbb{C}$-algebra maps,
and if two $\mathbb{C}^r$-algebra maps are conjugate as
$\mathbb{C}^r$-algebra maps then they are also conjugate as
$\mathbb{C}$-algebra maps. Hence it is suffient to prove that if $f:
B \to R$ is a $\mathbb{C}$-algebra map preserving the augmentations
then it is conjugate to a map of $\mathbb{C}^r$-algebras. Let $1_i$
denote the $i$th direct summand of the identity (i.e. the identity
arrow at the $i$th object) in either $B$ or $R$. Since $f$ preserves
the augmentations we must have
$$ f(1_i) = 1_i + m_i $$
for some elements $m_i \in \mathfrak{m}$, and since $f$ is an
algebra map we must have $\set{1_i + m_i}$ a complete orthogonal set
of idempotents. If we conjugate by the element
$$\sum_i 1_i+ m_i 1_i = 1 + \sum_i m_i 1_i \in R$$
(which is invertible since $\mathfrak{m}$ is nilpotent) then we get
a map that sends $1_i\in B$ to $1_i\in R$ for all $i$ and hence is a
map of $\mathbb{C}^r$-algebras.
\end{proof}

Now we consider the equivalent description (Corollary
\ref{Def_McontrolledbyExt}) of our deformation functor
$\mathcal{D}ef_\mathcal{M}$ as being the deformation functor
$\mathcal{D}ef_{Ext(\mathcal{M})}$ associated to the
$A_\infty$-category $\Ext_A(\mathcal{M})$.

Firstly recall the construction (discussed in the Introduction) of
an algebra from the degree one and two parts of an $A_\infty$
algebra. Let the higher products on $\Ext_A(\mathcal{M})$ be denoted
$m_i$. The degrees of the $m_i$ dictate that for all $i$
$$m_i : \Ext^1_A(\mathcal{M})^{\otimes i} \to \Ext^2_A(\mathcal{M})$$
The direct sum of all these maps is the homotopy Maurer-Cartan
function
$$HMC = \bigoplus_{i>0} m_i : T(\Ext^1_A(\mathcal{M})) \to \Ext^2_A(\mathcal{M})$$
The $m_i$ and $HMC$ are all maps of $\mathbb{C}^r$-bimodules. Now
given a $\mathbb{C}^r$-bimodule $V$, its $\mathbb{C}$-linear dual
$V^\vee=\Hom_\mathbb{C}(V,\mathbb{C})$ is also a
$\mathbb{C}^r$-bimodule. Dualizing the map $HMC$ in this way we get
a map
$$ HMC^\vee: \Ext^2_A(\mathcal{M})^\vee \to \hat{T}(\Ext^1_A(\mathcal{M})^\vee)$$
Quotienting by the two-sided ideal generated by the image gives us a
$\mathbb{C}^r$-algebra
$$\frac{\hat{T}(\Ext^1_A(\mathcal{M})^\vee)}{(\Ext^2_A(\mathcal{M})^\vee)}$$
This is naturally augmented  because there is no $m_0$ term in
$HMC$.

We claim that this algebra is (nearly) the ring of functions on the
formal deformation space associated to $\Ext_A(\mathcal{M})$. Our
use of the word \quotes{space} here is a little shaky, since this is
a non-commutative $\mathbb{C}^r$-algebra and we are not proposing to
define Spec of it! Nevertheless the statement makes rigourous sense
if we interpret it at the level of deformation functors.

At this point we must insert an extra assumption.

\begin{prop} \label{defMrepbyE} Assume that $\Ext^1_A(\mathcal{O}_{p_i}, \mathcal{O}_{p_j})$ is finite-dimensional for each $i,j$. Then for any $(R,\mathfrak{m})\in
\mathbf{Art}_\mathbb{C}^r$ there is a functorial isomorphism
$$\mathcal{D}ef_{Ext(\mathcal{M})}(R) = \mathbf{Alg}_\mathbb{C}^{\star r} \left(\frac{\hat{T}
(\Ext^1_A(\mathcal{M})^\vee)}{(\Ext^2_A(\mathcal{M})^\vee)},
R\right)\;/\; \set{\mbox{Inner $\mathbb{C}$-algebra automorphisms in
} R}$$
\end{prop}
\begin{proof}
Ignoring the gauge group for a minute, we have that
$\mathcal{D}ef_{Ext(\mathcal{M})}(R)$ is the zero locus of the
homotopy Maurer-Cartan function in
\begin{eqnarray}\label{needfindim}\Ext^1(\mathcal{M})\underline{\otimes} \mathfrak{m} &=& \bigoplus_{i,j}
\Hom_\mathbb{C} (\Ext^1_A(\mathcal{O}_{p_i},
\mathcal{O}_{p_j})^\vee, \mathfrak{m}_{ij})
\\ \label{needfindim2}
&=& \mathbf{Alg}_\mathbb{C}^{\star r} (
\hat{T}(\Ext^1_A(\mathcal{M})^\vee), R)\end{eqnarray} Here we have
used our finite-dimensionality assumption on the $\Ext^1$s. A map in
$\mathbf{Alg}_\mathbb{C}^{\star r} (
\hat{T}(\Ext^1_A(\mathcal{M})^\vee), R)$ is a zero of $HMC$
precisely when it induces a map
$$ \frac{\hat{T}
(\Ext^1_A(\mathcal{M})^\vee)}{(\Ext^2_A(\mathcal{M})^\vee)} \to R $$
on the quotient algebra.

The gauge group is the exponential of
$$ \Ext^0_A(\mathcal{M})\underline{\otimes} \mathfrak{m} = \bigoplus_i \mathfrak{m}_{ii}$$
which is $1 + \bigoplus_i \mathfrak{m}_{ii} \subset R $, and it acts
by conjugacy on $\mathbf{Alg}_\mathbb{C}^{\star r} (
\hat{T}(\Ext^1_A(\mathcal{M})^\vee), R)$. Now an arbitrary inner
$\mathbb{C}$-algebra automorphism of $R$ is a conjugation by an
element in
$$ (\mathbb{C}^*)^r + \bigoplus_i \mathfrak{m}_{ii}$$
but $(\mathbb{C}^*)^r$ is in the centre of $R$ so the orbits under
the gauge group are precisely
inner-$\mathbb{C}$-algebra-automorphism classes.

Functoriality follows from the functoriality of equations
(\ref{needfindim}) and (\ref{needfindim2}).
\end{proof}

We have nearly proved our main theorem. If we forget about inner
automorphisms, we have shown than $\mathcal{D}ef_\mathcal{M}$ is the
same as $\mathbf{Alg}_\mathbb{C}^{\star r}(A, -)$. This functor is
\textit{pro-representable}, i.e. it is represented by the completion
$\hat{A}_p$ of $A$ at the augmentation ideal (the kernel of $p$),
which can be thought of as a directed system of objects in
$\mathbf{Art}_\mathbb{C}^r$. On the other hand we know
$\mathcal{D}ef_\mathcal{M}=\mathcal{D}ef_{\Ext(\mathcal{M})}$, which
we have just shown to be pro-represented by
$$\frac{\hat{T}(\Ext^1_A(\mathcal{M})^\vee)}{(\Ext^2_A(\mathcal{M})^\vee)}$$
A pro-representing object is unique, so these two formal
$\mathbb{C}^r$-algbras must be the same.

Now we just have to check that this argument still holds when we
remember about the inner automorphisms, but this is just a matter of
carefully checking the standard proof of uniqueness of a
pro-representing object.

\begin{thm} \label{hatE=hatA}
Let $$\hat{E}
=\frac{\hat{T}(\Ext^1_A(\mathcal{M})^\vee)}{(\Ext^2_A(\mathcal{M})^\vee)}$$
then $\hat{E}$ is isomorphic as an augmented $\mathbb{C}^r$-algebra
to the completion $\hat{A}_p$ of $A$ at the kernel of $p$.
\end{thm}
\begin{proof}
We have shown (Corollary \ref{Def_McontrolledbyExt}, Proposition
\ref{defMrepbyA} and Proposition \ref{defMrepbyE}) that there is a
natural isomorphism $\Psi$ between
$$\frac{\mathbf{Alg}_\mathbb{C}^{\star r}(A, - )}{\mbox{Inner Aut.}} = \frac{\mathbf{Alg}_\mathbb{C}^{\star r}(\hat{A}_p, -
)}{\mbox{Inner Aut.}} : \mathbf{Art}_\mathbb{C}^r \to \mathbf{Set}
$$ and
$$\frac{\mathbf{Alg}_\mathbb{C}^{\star r}(\hat{E}, -
)}{\mbox{Inner Aut.}} : \mathbf{Art}_\mathbb{C}^r \to \mathbf{Set}
$$ Let $I \subset A$ be the kernel of $p$. Let $\pi_i: \hat{A}_p \to
A/I^i$ be the maps in the limiting cone on the diagram
\begin{equation} \label{Adiag} ...\onto A/I^3 \onto A/I^2 \onto A/I \end{equation}
Applying $\Psi$ to the cone $\set{\pi_i}$ (and picking
representatives of the resulting conjugacy classes) gives a set of
maps $\Psi\pi_i: \hat{E}\to A/I^i$ forming a cone that commutes up
to conjugacy. In fact since every map in (\ref{Adiag}) is a
surjection we may inductively pick representatives such that
$\Psi\pi_i$ forms a genuinely commuting cone. This cone then factors
through some map $f: \hat{E} \to \hat{A}_p$, and by naturality $\Psi
= \circ f$. Similarly, since $\hat{E}$ is also a limit of such a
diagram, there is a map $g: \hat{A}_p \to \hat{E}$ such that
$\Psi^{-1} = \circ g$.

The composition $\circ fg$ is the identity transformation, so
applying it to the cone $\set{\pi_i}$ we see that for each $i$ there
is an inner automorphism $\alpha_i$ of $A/I^i$ such that $\alpha_i
\pi_i fg = \pi_i$. We know $\pi_i$ is a surjection, so $\pi_i fg$
must also be a surjection, so by a quick diagram chase the maps
$\set{\alpha_i}$ commute with the maps in (\ref{Adiag}) and thus
lift to an automorphism $\tilde{\alpha}$ of $\hat{A}_p$. Then $\pi_i
\tilde{\alpha} fg = \pi_i$ for all $i$ and hence $\tilde{\alpha} fg
= \id_{\hat{A}_p}$. Similarly there is an automorphism
$\tilde{\epsilon}$ of $\hat{E}$ such that $\tilde{\epsilon} gf =
\id_{\hat{E}}$ so $f$ and $g$ must be isomorphisms.
\end{proof}

It has been suggested to the author by Lieven Le Bruyn (and
independently by Tom Bridgeland) that this result should generalise
to the case that $\mathcal{M}$ is a set of simple (not just
1-dimensional) modules, if we weaken isomorphism to Morita
equivalence.

\subsection{The graded case}\label{sectthegradedcase}

Let $A=A_\bullet$ be an $\mathbb{N}$-graded $\mathbb{C}^r$-algebra with $A_0=\mathbb{C}^r$.
The positively-graded part $A_{>0}$ of $A$ gives an augmentation, so we may
consider $A$ to be an object of $\mathbf{Alg}_\mathbb{C}^{\star r}$. We also consider $A_0$ as a set of $r$ one-dimensional $A$-modules.

Let $\hat{A}$ be the completion of $A$ at $A_{>0}$. The grading on
$A$, viewed as a $\mathbb{C}^*$ action, induces a $\mathbb{C}^*$ action on $\hat{A}$,
and we can recover $A \subset \hat{A}$ as the direct sum of the
eigenspaces of this action. Geometrically we may think of $A_0$ as a
repulsive fixed point (a source) for a $\mathbb{C}^*$ action on $A$, so if we take
an infinitesimal neighbourhood of the fixed point we can flow it outwards
until we see the whole of $A$.

There is an induced grading on $B(A, A_0)$ which we shall call a
\textit{lower grading} to distinguish it from usual dg structure (so
$B(A, A_0)$ is now bi-graded) and the differential obviously
preserves this lower grading. Dualizing we get the dg-category
$\mathcal{H} = \Hom_A(B(A, A_0), A_0)$ (see Section
\ref{sectdefthyofsetsofmods}). Now
$$B(A, A_0)^{-i} = A^{\otimes i+1}\otimes A_0 $$
so
$$\mathcal{H}^i = \Hom_\mathbb{C}( A^{\otimes i} , \mathbb{C})$$
so the grading on $A$ induces a splitting of $\mathcal{H}$ as a direct
product (not a direct sum) of lower graded pieces. The multiplication is degree zero
with respect to the lower grading.

\begin{lem} \label{Extisgraded} The category $\Ext_A(A_0)$ has an induced lower
grading, and there is a choice of $A_\infty$-structure on it such
that all the products, and the quasi-isomorphism $\Ext_A(A_0) \to
\mathcal{H}$, preserve the lower grading.
\end{lem}
\begin{proof}
$\Ext_A(A_0)$ aquires a lower grading since it is the homology of
$\mathcal{H}$ and the differential on $\mathcal{H}$ has lower degree
zero. Now we just apply the explicit form of the homological
perturbation algorithm (see e.g. \cite{merkulov}, \cite{lazaroiu}),
noting that since the differential and multiplication on
$\mathcal{H}$ have lower degree zero everything in the algorithm can
be chosen to respect the lower grading.
\end{proof}

The following theorem was proven in \cite{LPWZ} for the case $A_0 =
\mathbb{C}$. There they assume that $A$ is degree-wise finite-dimensional,
whereas we assume that $\Ext^1_A(A_0)$ is degree-wise
finite-dimensional. It is a consequence of the theorem that the two
assumptions are equivalent.

\begin{thm}\label{gradedresult}
Choose the $A_\infty$-structure on $\Ext_A(A_0)$ to be lower graded
as in the previous lemma, and assume $\Ext^1_A(A_0)$ is finite
dimensional in each lower degree. Then
$$ A = \frac{T(\Ext^1_A(A_0)^{\vee g})}{(\Ext^2_A(A_0)^{\vee g})}$$
as graded $\mathbb{C}^r$-algebras, where $\vee g$ denotes the lower-graded
$\mathbb{C}$-linear dual.
\end{thm}
\begin{proof}
First we note that the lower-grading ensures that
$$(HMC)^\vee : \Ext^2_A(A_0)^{\vee g} \to T(\Ext^1_A(A_0)^{\vee
g}) $$ so our statement makes sense. Without the grading it is
possible that the image of $HMC^\vee$ only lies in the completed
tensor algebra, as in Section \ref{sectdefingasetofpoints}.

Now we would like to use Theorem \ref{hatE=hatA}, but in Section
\ref{sectdefingasetofpoints} we required that $\Ext^1$ have finite
dimension, whereas now we are asking only for lower-degree-wise
finite-dimensionality. Examining the proofs however we see that as
long as we read $\vee g$ instead of $\vee$ then line
(\ref{needfindim}) still holds and hence Theorem \ref{hatE=hatA}
still holds. Thus if we complete both sides at their positively
graded parts then we have an isomorphism. This isomorphism respects
the lower grading however, since the isomorphism
$$ \mathcal{D}ef_{\Ext(A_0)} \cong \mathcal{D}ef_\mathcal{H} $$
respects the lower grading by Lemma \ref{Extisgraded} and the
isomorphism
$$\mathcal{D}ef_\mathcal{H} \cong \mathcal{D}ef_{A_0} $$
(Lemma \ref{defMcontrolledbyH}) clearly respects the lower grading.
Hence we can identify the original algebras on both sides since they
are the direct sums of the graded pieces of their completions.
\end{proof}

\section{Superpotential algebras} \label{sectsuperpotentials}

As discussed in the introduction, it is well known in the physics
literature that the algebras arising from quiver gauge theories on
Calabi-Yau three-folds can be described by a
\quotes{superpotential}. In this section (which probably has some
overlap with \cite{ginzburg}) we show why this is a consequence of
applying our results on deformation theory to the special case of
Calabi-Yau 3-folds.

Let $X$ be any complex manifold, and let $E$ be a vector bundle on
$X$ with holomorphic structure given by $\bar{\partial}$. Then all
other holomorphic structures on $E$ are given by adding to
$\bar{\partial}$ an element $$a \in \End(E)\otimes
\mathcal{A}_X^{0,1}$$ satisfying the Maurer-Cartan equation
$$\bar{\partial}(a) + a\wedge a = 0$$
and two such $a$ give isomorphic holomorphic structures if they
differ by a gauge transformation. This is dga deformation theory
again (see Section \ref{secttheAinftydeftheoryofapt}) - the
deformations of $(E,\bar{\partial})$ are governed by the dga
$$\RHom(E,E)\simeq \End(E,E)\otimes \mathcal{A}_X^{0,\bullet}$$
However, rather remarkably in this case the dga gives us the whole
moduli space, not just a formal neighbourhood.

When $X$ is a Calabi-Yau three-fold the Maurer-Cartan equation can
be written as the derivative of a (locally-defined) function: the
\textit{Chern-Simons} function
$$CS(a) := \int_X \Tr(\frac 12 a\wedge \bar{\partial}(a) + \frac 13 a\wedge a\wedge a
) \wedge \omega_{vol} $$ where $\omega_{vol}$ is a choice of
holomorphic volume form on $X$. Thus heuristically the moduli space
is the critical locus of this function. This means that we expect it
to be zero-dimensional, and that the number of points in it is the
Euler characteristic of the ambient space. Of course this is only
heuristic, since the ambient space is the quotient of an
infinite-dimensional vector space by an infinite-dimensional group.
What one can do however is to construct the moduli space using
algebraic geometry and then use the technology of $\textit{symmetric
obstruction theories}$, this leads to the definition of
Donaldson-Thomas invariants \cite{thomas}.

If we only care about formal deformations then we have a
finite-dimensional version of the above. As in Section
\ref{sectsketch}, we can replace the dga \linebreak$\End(E,E)\otimes
\mathcal{A}_X^{0,\bullet}$ by its homology equipped with an
$A_\infty$-structure, and the Maurer-Cartan equation by the Homotopy
Maurer-Cartan equation. This is still (formally) the critical locus
of a function - we just have to add all the higher products into the
Chern-Simons function.

Suppose now we have a 3-dimensional Calabi-Yau algebra $A$ instead
of a space. The same argument applies, so if $M$ is an $A$-module
then a formal neighbourhood of the moduli space of $M$ is the
critical locus of a function. In particular, if $M$ is the module
corresponding to a point of $A$ then a formal neighbourhood of $M$
in $A$ is the critical locus of a function. This function is called
a \textit{superpotential} for $A$ (at that point). We give a formal
definition of superpotentials in Section
\ref{sectsuperpotentialsdefinition}, and then make this argument
rigorous in Section \ref{sectthe3dimCYcase}.

In fact it is easy to construct global moduli spaces of modules over
$A$, and the construction is entirely finite-dimensional.
Furthermore if we have a (polynomial) superpotential for $A$, then
every moduli space of $A$-modules is globally the critical locus of
a function induced by the superpotential. This means that (rather
trivially!) the moduli spaces carry symmetric obstruction theories
and hence we can define invariants analogous to Donaldson-Thomas
invariants. We do this in Section \ref{sectDT}.

In the examples discussed in the Introduction we have both a space
$X$ and an algebra $A$, and they are derived equivalent, so there is
some relationship between moduli spaces of sheaves on $X$ and moduli
spaces of $A$-modules. The physical picture (as we discussed in
Section \ref{sectphysicalargument}) is that moving from $X$ to $A$
means that we are moving in the stringy K\"ahler moduli space, which
mathematically probably means some space of stability conditions on
the triangulated category $D^b(X)$. This space (or at least a
conjectural version of it) has been constructed by Bridgeland
\cite{bridgeland2}. Part of this picture is obvious: passing from
$D^b(X)$ to the equivalent $D^b(A)$ is just a change of T-structure,
which is part of a Bridgeland stability condition. The remaining
data, called the central charge, should roughly correspond on the
algebra side to putting a GIT stability condition on the moduli
space of $A$-modules. It should be possible to build a function over
the whole of the space of stability conditions such that you can
Taylor expand it at the point corresponding to $X$ or the point
corresponding to $A$ and get the generating function of the
Donaldson-Thomas invariants of $X$ or the Donaldson-Thomas-type
invariants of $A$ respectively. This has been carried out in one
example by Szendr\"oi \cite{szendroi}, and much interesting work is
being done (e.g. \cite{joyce}, \cite{cs}), and much remains to be
done, to properly understand this picture.

\subsection{Definition of superpotentials}\label{sectsuperpotentialsdefinition}
{\large } Let $V$ be a $\mathbb{C}^r$-bimodule. A
\textit{superpotential} is simply a sum of cycles in the path
algebra of $V$ taken up to cyclic permutation, i.e. an element
$$W \in TV/\;[TV, TV]$$
If we instead use the completed tensor algebra we have a
\textit{formal superpotential}, i.e. an element of $\hat{T}V /\;
[\hat{T}V, \hat{T}V]$. Note that in the introduction we identified
these spaces with the spaces of cyclicly symmetric elements in $TV$
and $\hat{T}V$.

Roughly, we are interested in the (non-commutative) affine scheme
described by the critical locus of $W$, but we have to take a little
care with our definition of partial derivative.

\begin{defn} \label{cycpartderiv}Let $x\in V^\vee$. For any $t$ we define the \textit{cyclic partial derivative}
in the direction of $x$ as the map
$$\partial^\circ_x : V^{\otimes t} \to V^{\otimes t-1} $$
$$\partial^\circ_x(v_1\otimes...\otimes v_t) = \sum_{s=1}^t x(v_s)
v_{s+1}\otimes...\otimes v_t\otimes v_1\otimes ... \otimes v_{s-1}
$$
Taking direct sums/products we get maps
$$ \partial^\circ_x: TV/\;[TV, TV] \to TV$$
and
$$\partial^\circ_x: \hat{T}V/\;[\hat{T}V,\hat{T}V] \to \hat{T}V$$
\end{defn}

\begin{defn} The algebra \textit{generated by a superpotential} $W$
is $TV/(R)$ where $R$ is the subspace
$$ R = \partial^\circ W := \set{ \partial^\circ_x W \;\;|\;\; x\in V^\vee } \subset TV$$
If $W$ is a formal superpotential then it generates the algebra
$\hat{T}V/(R)$ with the same definition of $R$.
\end{defn}

\subsection{3-dimensional Calabi-Yau algebras}
\label{sectthe3dimCYcase}

As in Section \ref{sectdefingasetofpoints}, we pick an augmented
$\mathbb{C}^r$-algebra $(A,p) \in \mathbf{Alg}_\mathbb{C}^{\star r}$  and assume that the
resulting set of $r$ one-dimensional $A$-modules
$\mathcal{M}=\set{\mathcal{O}_{p_i}}$ has
$\Ext^1_A(\mathcal{O}_{p_i}, \mathcal{O}_{p_j})$ finite dimensional
$\forall i,j$.

Recall that an $A_\infty$ category $(\mathcal{C},m_i)$ is
\textit{Calabi-Yau of dimension} $d$ if it carries a trace map
$$\Tr_M: \Hom(M,M) \to \mathbb{C}$$
of degree $-d$ for each object $M\in \mathcal{C}$ which is closed
with respect to the differential on $\Hom(M,M)$, such that the
multilinear maps
$$ \Tr_M\circ m_i : \Hom(M, M_1)\otimes\Hom(M_1,M_2)\otimes...\otimes
\Hom(M_{i-1}, M) \to \mathbb{C} $$ are (graded) cyclically
symmetric. The bilinear pairing
$$\langle\;\;\rangle_{M,N} := \Tr_M\circ m_2: \Hom(M,N)\otimes \Hom(N,M) \to \mathbb{C}$$
is required to be non-degenerate on homology.

The canonical example is the derived category $D^b(X)$ of a compact
smooth Calabi-Yau variety $X$ with the Serre duality pairing (see
e.g. \cite{costello2}, Section 2.2). We take a dg-model of the
derived category whose objects are finite complexes of vector
bundles, with morphisms between complexes $E^\bullet$ and
$F^\bullet$ given by the complex
$$ \Hom(E^\bullet, F^\bullet \otimes \Omega_X^{0,*})$$
where $\Omega_X^{0,*}$ is the complex of $(0,*)$-forms with the
exterior differential. The Calabi-Yau pairing is
$$\langle \alpha | \beta \rangle  = \int_X
\Tr(\alpha\wedge\beta)\wedge vol $$ where $vol$ is a holomorphic
volume form. Cyclic symmetry follows from the symmetry of $\Tr$ and
$\wedge$. It is then possible to run the homological perturbation
algorithm for this example in such a way that cyclic symmetry is
preserved (\cite{lazaroiu}), so we get a Calabi-Yau
$A_\infty$-structure on $D^b(X)$.

 For our non-compact
examples in Section \ref{sectlocalCYs} we will get cyclicity by
another construction. In general the cyclic symmetry of a pairing on
an $A_\infty$-category seems to be a delicate notion (e.g. it is not
clear that it is preserved by homological perturbation), and should
maybe be relaxed to some homotopy invariant notion (\cite{KS},
Chapter 10).

\begin{thm}\label{CY3impliesSP} Suppose that $\langle \mathcal{M}\rangle
\subset D^b(A)$ is Calabi-Yau of dimension 3. Then the completion
$\hat{A}_p$ of $A$ at the kernel of $p$ is given by a formal
superpotential.
\end{thm}

\begin{proof}
We learnt this construction from \cite{lazaroiu}. Consider the
formal superpotential
 $$W \in \hat{T}(\Ext^1_A(\mathcal{M})^\vee) \;/\;
 \; [\hat{T}(\Ext^1_A(\mathcal{M})^\vee), \hat{T}(\Ext^1_A(\mathcal{M})^\vee)] $$
defined by
\begin{equation}\label{superpotential}
W(a_1\otimes...\otimes a_t) := \frac{1}{t}\left\langle\; a_1 \;|\; m_{t-1}(a_2\otimes
...\otimes a_t)\; \right\rangle
\end{equation}
where  $\set{m_i}$ are the $A_\infty$ products on
$\Ext_A(\mathcal{M})$ and $\langle\cdot | \cdot\rangle$ denotes the
Calabi-Yau pairing on $\langle\mathcal{M}\rangle$. This is
well-defined by the cyclicity of $m_t$. Then for $x\in
\Ext^1_A(\mathcal{M})$ we have
\begin{eqnarray*}
\partial^\circ_x W (a_1\otimes...\otimes a_t) & =& \sum_{s=1}^{t+1} W(a_s\otimes...
\otimes a_t\otimes x\otimes a_1\otimes...\otimes a_{s-1})\\
&=& \langle \; x \;|\; m_t(a_1\otimes...\otimes a_t)\; \rangle
\end{eqnarray*}
However, if we recall the definition of the $HMC$ function
$$HMC = \bigoplus_{i>0} m_i : T(\Ext^1_A(\mathcal{M})) \to \Ext^2_A(\mathcal{M})
$$
we see that we have shown
$$\partial^\circ_x W = HMC^\vee(x)$$
under the identification $\Ext^2_A(\mathcal{M})^\vee = \Ext^1_A(\mathcal{M})$
given by $\langle\cdot | \cdot\rangle$. Hence by Theorem
\ref{hatE=hatA} the algebra generated by $W$ is $\hat{A}_p$.\end{proof}

Combining this construction with what we know about graded algebras
(Theorem \ref{gradedresult}) we recover a result of Bocklandt.

\begin{thm} \label{bocklandt}\cite{bocklandt} Let $A$ be a graded degree-wise finite-dimensional
algebra with $A_0$ semi-simple and with the subcategory of $D^b(A)$
generated by the summands of $A_0$ Calabi-Yau of dimension 3. Then
$A$ is given by a superpotential.
\end{thm}

Note that we do not require that the whole of $D^b(A)$ be
Calabi-Yau. Indeed this will not be the case for the algebras coming
from non-compact Calabi-Yau varieties that we study in Section
\ref{sectlocalCYs} below. Also note that the proof in
\cite{bocklandt} constructs a cyclic $A_\infty$ structure directly
rather than assuming it.

\subsection{Donaldson-Thomas-type invariants}\label{sectDT}

Let $A$ be a $\mathbb{C}^r$-algebra with generators $V$ and
relations $R \subset TV$. The moduli space of $A$-modules and the
role of GIT stability conditions for it was explained in
\cite{king}, the following is a summary.

The \textit{dimension vector} of an $A$-module $M$ is just the
vector
$$\mathbf{d} = (d_1,...,d_r)$$
where $d_i$ is the dimension of the $i$th summand $M_i$ of $M$. To
give an $A$-module of dimension $\mathbf{d}$ we have to give a
representation of $A$ on the $\mathbb{C}^r$-module
$$\mathbb{C}^\mathbf{d} = \bigoplus_i \mathbb{C}^{d_i}$$
To start with forget about the relations, and consider the set of
representations of $TV$ on $\mathbb{C}^\mathbf{d}$. These form the
$\mathbb{C}^r$-bimodule
$$V^\vee \underline{\otimes} \End_\mathbb{C} (\mathbb{C}^\mathbf{d})$$
since an element of this space is precisely a linear map
$$V_{ij} \to \Hom_\mathbb{C}(\mathbb{C}^{d_i}, \mathbb{C}^{d_j})$$
for all $i,j$. Now a relation $r\in R$ is an element of $TV$, and so
it induces, using the composition in $\End_\mathbb{C} (\mathbb{C}^\mathbf{d})$, a
polynomial map
$$ r: V^\vee \underline{\otimes} \End_\mathbb{C} (\mathbb{C}^\mathbf{d}) \to \End_\mathbb{C} (\mathbb{C}^\mathbf{d}) $$
The zero locus of $r$ is just those representations that obey the
relation $r$. An $A$-module structure on $\mathbb{C}^\mathbf{d}$ is a
representation that obeys all the relations, so the scheme of
$A$-module structures on $\mathbb{C}^\mathbf{d}$ is the common zero locus $Z$
of all such $r\in R$.

 Two such modules are isomorphic if they differ
by a change of basis in $\mathbb{C}^d$, which is an element of the group
$$GL(\mathbf{d}):= \prod_i GL(d_i, \mathbb{C})$$
Hence the moduli stack of $A$-modules of dimension $\mathbf{d}$ is
$$\mathcal{M}_{A,\mathbf{d}} = [\,Z\, /\, GL(\mathbf{d})\,] $$
If we pick a character $\xi$ of $GL(\mathbf{d})$ then we can
instead take the GIT quotient
$$\mathcal{M}_{A,\mathbf{d}}^\xi = Z \sslash{\xi} GL(\mathbf{d}) $$
A character of $GL(\mathbf{d})$ is necessarily of the form
$$\xi(g) = \prod_{i=1}^r \det(g_i)^{\theta_i}$$
for some $r$-tuple of integers $(\theta_i)$. Hence given a character
$\xi$ we can define a function
$$\Theta_\xi: K_0(A\mathbf{-mod}) \to \mathbb{Z}$$
by sending a module of dimension vector $\mathbf{d}$ to the integer
$$\sum_i \theta_i d_i $$
This function is indeed well defined on $K_0$ because each $d_i$ is
additive over short exact sequences.

\begin{defn}\cite{king} Let $\Theta: K_0(A\mathbf{-mod}) \to \mathbb{R}$ be an additive
function. A module $M$ is $\Theta$-semistable (resp.
$\Theta$-stable) if $\Theta(M) = 0$ and every proper submodule
$N\subset M$ has $\Theta(N)\geq 0$ (resp. $>0$).\end{defn}

Two $\Theta$-semistable modules are called \textit{S-equivalent} if
they have the same composition factors in the abelian category of
$\Theta$-semistable modules.

\begin{thm}\cite{king} $\mathcal{M}_{A,\mathbf{d}}^\xi$ is a coarse moduli
space for $\Theta_\xi$-semistable $A$-modules up to S-equivalence.
\end{thm}

We say $\mathbf{d}$ is \textit{indivisible} if it is not a multiple
of another integral vector.

\begin{thm}\cite{king}\label{M_Aisfinemodspace} If $\mathbf{d}$ is indivisible and there are no
strictly $\Theta_\xi$-semistable modules of dimension $\mathbf{d}$
then $\mathcal{M}_{A,\mathbf{d}}^\xi$ is a fine moduli space for
$\Theta_\xi$-stable $A$-modules.
\end{thm}

Now we consider this construction when $A$ is a superpotential
algebra given by a superpotential $W$. We will show that in this
case $\mathcal{M}_{A,\mathbf{d}}^\xi$ carries a natural symmetric
obstruction theory.

It is common in algebraic geometry that moduli spaces fail to have
the \quotes{expected} dimension that a naive calculation predicts.
If one then tries to produce invariants by integrating cochains over
the moduli space then one gets unhelpful results, because the
cochains of the \quotes{correct} degree are not top-degree and
integrate to zero. Obstruction theories (introduced by \cite{BF} and
others) are pieces of technology that resolve this problem (at least
if the actual dimension is bigger than the expected dimension), by
producing a \textit{virtual fundamental class} on the moduli space,
of degree equal to the expected dimension. One can then integrate
against this, instead of against the usual fundamental class.

The intuition behind obstruction theories is simple. Suppose that
your moduli space $X$ is cut out of an ambient space $Y$ by some
section $\sigma\in \Gamma(E)$ of a vector bundle. The expected
dimension of $X$ is then dim$(Y)-$rank$(E)$, and if $\sigma$ is
transverse then this is also the actual dimension. Suppose now that
$\sigma$ is not transverse, then the derivative of $\sigma$ at its
zero locus $X$ gives an exact sequence
$$0 \to TX \to TY|_X \stackrel{D\sigma}{\longrightarrow} E|_X \to Obs \to 0 $$
where $Obs$ (the \textit{obstruction sheaf}) is some sheaf on $X$
given by the cokernel of $D\sigma$. If $Obs$ is a vector bundle then
we can imagine perturbing $\sigma$ by adding on a small transverse
section $\tau\in\Gamma(Obs)$. Then $\sigma+\tau$ is transverse, and
its zero locus is the zero locus of $\tau$ in $X$, which is the
euler class of $Obs$. This has degree equal to the expected
dimension. Now we can abstract this: an obstruction theory on $X$ is
(roughly) a sheaf $Obs$ on $X$ and an exact sequence
$$ 0 \to TX \to E_1 \to E_2 \to Obs \to 0$$
where $E_1$ and $E_2$ are vector bundles. The associated virtual
fundamental class is $eu(Obs)$ if $X$ is smooth (see \cite{BF} for
the general technology).

A \textit{symmetric} obstruction theory corresponds to the special
case when $E = T^{*}X$ and $\sigma$ is the derivative of a function.
This means that the exact sequence given by $D\sigma$ is self-dual,
since the Hessian of a function is a symmetric matrix. Symmetric
obstruction theories were introduced in $\cite{behrend}$, and are
the right technology for Donaldson-Thomas invariants (as explained
at the start of Section \ref{sectsuperpotentials}).

However, in our case much of this advanced technology is redundant,
since what we are going to show is that our moduli spaces
$\mathcal{M}_{A,\mathbf{d}}^\xi$ are genuinely the critical locus of
some function on a finite dimensional space.

 Recall that $W$ is required to be
a sum of cycles in $TV$. Using the composition in
$\End_\mathbb{C}(\mathbb{C}^\mathbf{d})$ it induces a map
$$W: V^\vee \underline{\otimes} \End_\mathbb{C} (\mathbb{C}^\mathbf{d}) \to \bigoplus_i \End_\mathbb{C}
(\mathbb{C}^{d_i})$$
Now we can take traces at each of the $r$ vertices and
sum them, getting a scalar polynomial function
$$ \tilde{W} : V^\vee \underline{\otimes} \End_\mathbb{C} (\mathbb{C}^\mathbf{d}) \to \mathbb{C}$$
$$ \tilde{W} = \Tr(W) $$
\begin{prop}\label{ZiscritlocofW}The (scheme-theoretic) critical locus of $\tilde{W}$ is precisely the
(scheme-theoretic) zero locus $Z$ of the relations.
\end{prop}

This is just the statement that the partial derivatives of $\tilde{W}$ are
the polynomials on $V^\vee \underline{\otimes} \End_\mathbb{C} (\mathbb{C}^\mathbf{d})$ induced
by the relations.

\begin{proof}
Pick the standard basis of $\mathbb{C}^\mathbf{d}$ so that elements of
$\End_\mathbb{C} (\mathbb{C}^\mathbf{d})$ are matrices. The heart of the proof is just
the fact that if we take two independent matrices $M$ and $N$ and
then partially differentiate the function
$$\Tr(MN)$$
holding $N$ fixed, we get the matrix $N^T$. More generally if
$\set{M_1,...,M_l}$ are independent matrices and we partially
differentiate the function
\begin{equation}\label{TrM}
\Tr(M_{i_1}...M_{i_t}) \end{equation}
 by varying $M_j$ we get the
transpose of the matrix
\begin{equation}\label{derivofTrM}\sum_{i_s=j}
M_{i_{s+1}}...M_{i_t}M_{i_1}...M_{i_{s-1}} \end{equation}
 Now pick a
basis $\set{e_1,..,e_l}$ of $V$ and a let the dual basis of $V^\vee$
be $\set{\epsilon_1,...,\epsilon_l}$. Then an element $M\in V^\vee
\underline{\otimes} \End_\mathbb{C} (\mathbb{C}^\mathbf{d})$ is given by a set of
matrices $\set{M_1,...,M_l}$, and evaluating the function
$\tilde{W}$ at $M$ gives a linear combination of terms of the form
of (\ref{TrM}). Thus taking partial derivatives of $\tilde{W}$ in
all of the $M_j$ directions gives a function
$$V^\vee \underline{\otimes} \End_\mathbb{C} (\mathbb{C}^\mathbf{d}) \to  \End_\mathbb{C}
(\mathbb{C}^\mathbf{d}) $$ which is the transpose of the corresponding sum of
terms of the form (\ref{derivofTrM}). If we recall our Definition
\ref{cycpartderiv} of the cyclic partial derivative we see that this
function is the transpose of
$$\partial^\circ_{\epsilon_j} W : V^\vee \underline{\otimes} \End_\mathbb{C} (\mathbb{C}^\mathbf{d}) \to  \End_\mathbb{C}
(\mathbb{C}^\mathbf{d}) $$ Since the set $R$ of relations is spanned by
$\partial^\circ_{\epsilon_j} W$ the proposition is proved.
\end{proof}

\begin{cor} If there are no
strictly $\xi$-semistable points in $V^\vee \underline{\otimes} \End_\mathbb{C} (\mathbb{C}^\mathbf{d})$ then $\mathcal{M}^\xi_{A,\mathbf{d}}$ carries a symmetric
obstruction theory.
\end{cor}
\begin{proof}
Since there are no strictly semistables the ambient space
$$\mathcal{A}^{\xi}_\mathbf{d} := V^\vee \underline{\otimes} \End_\mathbb{C} (\mathbb{C}^\mathbf{d}) \sslash{\xi} GL(\mathbf{d})$$
is smooth. By invariance, $\tilde{W}$ descends to a function on $\mathcal{A}^\xi_\mathbf{d}$
 whose critical locus is $\mathcal{M}^\xi_{A,\mathbf{d}}$. The Hessian of this function gives a symmetric obstruction theory.
\end{proof}

 Associated to this obstruction theory is a virtual fundamental class
 $[\mathcal{M}^\xi_{A,\mathbf{d}}]^{vir}$.

\begin{defn} Let $\mathbf{d}$ be indivisible, and assume (i) there are no
strictly $\Theta_\xi$-semistable modules of dimension $\mathbf{d}$, and
(ii) $\mathcal{M}^\xi_{A,\mathbf{d}}$ is compact. Then we define the \textit{Donaldson-Thomas-type invariant}
$$\tilde{N}_{A, \mathbf{d}, \Theta_\xi} = \int [\mathcal{M}^\xi_{A,\mathbf{d}}]^{vir}$$
\end{defn}

In light of Theorem \ref{M_Aisfinemodspace} this really is an invariant of
the pair $(A, \Theta_\xi)$, it does not depend on our presentation of $A$
or even on our choice of $r$ idempotents (in fact it should probably be thought
of as an invariant of $(A\mathbf{-mod}, \Theta_\xi)$ but we do not know how to
make this precise). Furthermore, by the usual obstruction theory arguments it is invariant
under deformations of $W$ that leave $\mathcal{M}^\xi_{A,\mathbf{d}}$ compact.

It might appear that $\tilde{N}_{A, \mathbf{d}, \Theta_\xi}$ depends on the obstruction theory. However, as Behrend has shown in
\cite{behrend}, the virtual count under a \textit{symmetric}
obstruction theory is in fact an intrinsic invariant equal to a
weighted Euler characteristic
$$\chi(X,\nu_X) = \sum_{n\in\mathbb{Z}} n\,\chi(\set{\nu_X = n})$$
where $\nu$ is a constructible function defined by Behrend that
exists on any DM stack and measures the singularity of the space.
Using this we can drop the compactness assumption and define
$$\tilde{N}_{A, \mathbf{d}, \Theta_\xi} = \chi(\mathcal{M}^\xi_{A,\mathbf{d}},
\nu)$$ though if $\mathcal{M}^\xi_{A,\mathbf{d}}$ is not compact
then we should not expect this to be deformation invariant.

We hope to pursue these invariants further in future work. For the moment
however we content ourselves with the following observation.

\begin{lem} Assume there are no strictly $\xi$-semistable points in $V^\vee \underline{\otimes} \End_\mathbb{C} (\mathbb{C}^\mathbf{d})$. Suppose that
$\mathcal{M}^\xi_{A,\mathbf{d}}$ is confined (scheme-theoretically) to a
compact submanifold
$$\mathcal{N}\subset \mathcal{A}^\xi_\mathbf{d}$$
Then $\mathcal{M}^\xi_{A,\mathbf{d}}$ is smooth.
\end{lem}
\begin{proof}
$\tilde{W}$ is a holomorphic function on $\mathcal{A}^\xi_\mathbf{d}$ so
it is constant along $\mathcal{N}$, so $d\tilde{W}$ restricted
to $\mathcal{N}$ is a section of the conormal bundle $N^\vee_\mathcal{N}$.
The zero locus of this section is $\mathcal{M}^\xi_{A,\mathbf{d}}$. At any
point of $\mathcal{M}^\xi_{A,\mathbf{d}}$ the symmetric obstruction theory
(i.e. the Hessian of $\tilde{W}$) is an exact sequence
$$0 \to T\mathcal{M}^\xi_{A,\mathbf{d}} \to T\mathcal{N} \oplus N_\mathcal{N} \stackrel{Dd\tilde{W}}{\longrightarrow} T^\vee\mathcal{N}\oplus N^\vee_\mathcal{N} \to Ob_{\mathcal{M}^\xi_{A,\mathbf{d}}} \to 0$$
By assumption  $T\mathcal{M}^\xi_{A,\mathbf{d}} \subset T\mathcal{N}$, so
$N_\mathcal{N} \stackrel{Dd\tilde{W}}{\longrightarrow} T^\vee\mathcal{N}$ is an injection, and dually $T\mathcal{N}\stackrel{Dd\tilde{W}}{\longrightarrow} N^\vee_\mathcal{N}$
must be a surjection. Hence $d\tilde{W}|_\mathcal{N}$ is a transverse section
of the conormal bundle, and $\mathcal{M}^\xi_{A,\mathbf{d}}$ is smooth.
\end{proof}

\section{The derived categories of some local Calabi-Yaus}
\label{sectlocalCYs}

\subsection{Ext algebras on local Calabi-Yaus}

 For any smooth scheme $Z$, there is a \quotes{formal} way
to extend $Z$ to a Calabi-Yau, namely we use the embedding
$$\iota: Z \into \omega_Z$$
of $Z$ as the zero section in the total space of its canonical
bundle. In this section we prove that this procedure is reflected at
the level of ($A_\infty$-enriched) derived categories, i.e. that
$$\iota_* D^b(Z) \subset D^b(\omega_Z)$$
is a formal Calabi-Yau enlargement of $D^b(Z)$. There is a slight
subtlety here: $D^b(\omega_Z)$ is not actually Calabi-Yau since
$\omega_Z$ is non-compact, however if $Z$ is compact then objects in
$\iota_* D^b(Z)$ are compactly supported, so this subcategory
\textit{is} Calabi-Yau. In any case, let us first explain what we
mean by \quotes{formal Calabi-Yau enlargement}.

Let $m:V\otimes V \to V$ be any bilinear map on a vector space.
We claim that $m$ naturally extends to a bilinear map
$$m^c: (V\oplus V^\vee)^{\otimes 2} \to (V \oplus V^\vee)$$
such that the associated trilinear dual
$$\tilde{m}^c: (V\oplus V^\vee)^{\otimes 3} \to \mathbb{C}$$
is cyclically symmetric. This is straightforward: we let $m^c$ be the direct
sum of $m$ with the two maps
$$m_1:V\otimes V^\vee  \to V^\vee$$
and
$$m_2:  V^\vee\otimes V  \to V^\vee$$
obtained by dualising and cyclically permuting $m$ (on the fourth direct
summand we declare $m^c$ to be zero). We shall call $m^c$
the \textit{cyclic completion} of $m$. Of course this construction is
hardly profound, so no doubt it has been studied and named already,
but unfortunately we do not have a reference for it.

In the same way we may cyclicly complete any collection of $n$-linear maps. Furthermore
we claim that this process is sufficently natural that any algebraic structure
present in the set of maps will be preserved. For example, if
$$m:V\otimes V \to V$$
is an associative unital product, then one easily checks that $m^c$ is associative
and inherits the unit of $m$. In fact in this case our construction is nothing
more than the extension algebra associated to the $(V,m)$-bimodule $V^\vee$,
but the point is that for this particular bimodule the extension algebra
is a Frobenius algebra. The general statement is the following:

\begin{prop} Let $\mathcal{P}$ be a cyclic operad in $\mathbf{dgVect}$ or
$\mathbf{dg}$-$\mathbb{C}^r$-$\mathbf{bimod}$, and let $U^*\mathcal{P}$
be the underlying classical operad. Let
$$\Psi: U^*\mathcal{P} \to \mathcal{E}nd_V$$
be an algebra over $U^*\mathcal{P}$ whose underlying dg
vector-space is $V$. Then $V\oplus
V^\vee$ is naturally an algebra over $\mathcal{P}$.
\end{prop}
For the definitions of cyclic and classical operads see
\cite{costello1}.
\begin{proof}
We construct a natural transformation of classical operads
$$\Phi: U^*\mathcal{P} \to \mathcal{E}nd_{V\oplus V^\vee}$$
by intertwining the procedure above with the action of the cyclic groups
on $\mathcal{P}$ (for the operad of Frobenus algebras this action is trivial). Let the components of
$\Psi$ be
$$\Psi_n: U^*\mathcal{P}(n) \to \Hom(V^{\otimes n}, V)$$
and let $\gamma_n$ be a generator of the cyclic group $C_{n+1}$. We define
$$\Phi_n: U^*\mathcal{P}(n) \to \Hom((V\oplus V^\vee)^{\otimes n}, V\oplus
V^\vee)$$
to be the direct sum of $\Psi_n$ and all the compositions
$$U^*\mathcal{P}(n) \xrightarrow{\Psi_n \gamma_n^{-k}} \Hom(V^{\otimes n}, V)
\into \Hom(V^{\otimes k-1}\otimes V^\vee \otimes V^{\otimes n-k}, V^\vee)$$
This procedure commutes with the process of gluing maps together along graphs,
so it does define a natural transformation of operads.

Pick a pairing $\rho$ on $(V\oplus V^\vee)^\vee$ compatible with the
natural pairing on $V\oplus V^\vee$ (if $V$ has
finite-dimensional homology then $\rho$ is unique up to homotopy, so
this choice should not worry us). Then $\mathcal{E}nd_{V\oplus
V^\vee}^\rho$ is a cyclic operad and there is a natural
transformation
$$\mathcal{E}nd_{V\oplus V^\vee} \to U^*\mathcal{E}nd_{V\oplus V^\vee}^\rho$$
The composition of this with $\Phi$ is equivariant with respect to the actions
of the cyclic groups by construction, so it lifts to a map
$$ \mathcal{P} \to \mathcal{E}nd_{V\oplus V^\vee}^\rho$$
\end{proof}
The lemma remains true if we fix the dimension of the pairing by
forming the \textit{n-dimensional cyclic completion}
$$V\oplus V^\vee[-n]$$
In particular we may form the $n$-dimensional cyclic completion of
an $A_\infty$-algebra.

\begin{thm} \label{Extonomegaiscycliccompletion} Let $Z$ be a smooth proper scheme of dimension $n-1$, and let
$$\iota: Z \to \omega $$
be the embedding of $Z$ into its canonical bundle. Then for any $S
\in D^b(Z)$, the $A_\infty$-algebra
$$\Ext_\omega (\iota_*S, \iota_*S) $$
is the $n$-dimensional cyclic completion of $\Ext_Z (S,S)$.
\end{thm}

\begin{proof}
Let $\pi: \omega \to Z$ be the projection. We have a tautological
exact sequence (which we draw right to left to ensure a happy
typographical coincidence later on):
\begin{equation}\label{tautexactseq} 0 \from \iota_*\mathcal{O}_Z \from \mathcal{O}_\omega \xleftarrow{\tau} \pi^*\omega^\vee \from 0 \end{equation}
so there are quasi-isomorphisms of chain complexes
\begin{eqnarray}\label{cyccompl}
 \RHom_\omega(\iota_* S, \iota_* S) & =& \RHom_\omega (
\pi^*S \xleftarrow{\tau} \pi^*(S\otimes \omega^\vee),\;\; \iota_* S)
\\\label{cyccompl2} & = & \RHom_Z(S,S) \oplus  \RHom_Z(S\otimes
\omega^\vee, S)[-1]
\end{eqnarray}
since $\tau$ vanishes along the zero section. This latter admits a
dga structure if we identify $\RHom_Z(S,S)$ with $\RHom_Z(S\otimes
\omega^\vee, S\otimes \omega^\vee)$ and declare the product of two
elements in $\RHom(S\otimes \omega^\vee,S)$ to be zero. This dga
structure is cyclic with respect to the Serre duality pairing, so it
is in fact the $n$-dimensional cyclic completion of $\RHom_Z(S,S)$.

We claim that under this dga structure the equation (\ref{cyccompl2}) is actually a dga quasi-isomorphism. To see this take
a dga model of the LHS of the form
$$\RHom_\omega (
\pi^*S \xleftarrow{\tau} \pi^*(S\otimes \omega^\vee),\;\;\pi^*S \xleftarrow{\tau} \pi^*(S\otimes \omega^\vee))$$
What we mean here is that we should apply $\RHom$ termwise. This gives us
a two-by-two term complex that looks like the square on the RHS of the following
diagram (we suppress the $\pi^*$'s to keep the width manageable):
 \begin{equation*}\xymatrix{
  \RHom_Z(S,S)\ar[d]^0  & \RHom_\omega(S, S)\ar[d]^\tau \ar@{-->}[l] & \RHom_\omega(S,S\otimes \omega^\vee)[1]  \ar[l]_\tau \ar[d]^\tau   \\
  \RHom_Z(S\otimes\omega^\vee,S)[-1]& \RHom_\omega(S\otimes \omega^\vee,S)[-1]
  \ar@{-->}[l]&  \RHom_\omega(S\otimes\omega^\vee, S\otimes\omega^\vee) \ar[l]_{-\tau}
   }
\end{equation*}
The dga structure on this two-by-two term complex is easy to describe (by
a happy typographical coincidence...): treat an element of the RHS square in the above diagram as a two-by-two matrix, then the product is precisely matrix multiplication. The dashed arrows are
the cokernels of each row, we know what these are from the tautological exact
sequence (\ref{tautexactseq}). The two dashed arrows together form the quasi-isomorphism of (\ref{cyccompl2}).

The dashed arrows do not form a map of dgas. However, consider the map going
in the opposite direction, from the LHS two-term complex to the RHS two-by-two
term complex, obtained by summing the maps
$$\pi^*: \RHom_Z(S,S) \to \RHom_\omega(S, S) $$
 $$\pi^*: \RHom_Z(S,S) \to \RHom_\omega(S\otimes\omega^\vee, S\otimes\omega^\vee)$$
and
$$\pi^*: \RHom_Z(S\otimes \omega^\vee, S)[-1] \to \RHom_\omega(S\otimes \omega^\vee,
S)[-1]$$ It is easy to check that this map respects the product
structures on each side, so it is a map of dgas. It is also a right
inverse to the quasi-isomorphism given by the dashed arrows, hence
it is a dga quasi-isomorphism. We conclude that
$\RHom_\omega(\iota_*S, \iota_*S)$ is the $n$-dimensional cyclic
completion of $\RHom_Z(S,S)$.

Now we apply the homological perturbation algorithm (e.g.
\cite{merkulov}, \cite{lazaroiu}). We have to pick maps
\begin{equation*}\xymatrix{
  \Ext_\omega(\iota_*S, \iota_*S)\ar@<.5ex>[r]^i  & \RHom_\omega(\iota_*S, \iota_*S) \ar@<.5ex>[l]^p   \\
   }
\end{equation*}
and a homotopy $h$ such that $pi = \id$ and $ip = \id +
\partial(h)$. We can pick $i, p$ and $h$ such that they are the
cyclic completions of corresponding maps between $\Ext_Z(S,S)$ and
$\RHom_Z(S,S)$, it is then clear that running the algorithm on
$\RHom_\omega(\iota_*S, \iota_*S)$ is the same as running it on
$\RHom_Z(S,S)$ and then cyclically completing.
\end{proof}

\subsection{Completing the algebra of an exceptional collection}

Now we specialize to the case discussed in the introduction, where
$Z$ is a surface and $\omega$ is a local CY
three-fold. Furthermore we assume that we have been given a finite full
strong exceptional  collection of objects $\set{T_i} \subset D^b(Z)$. Known
examples include the del Pezzo \cite{KO} and ruled surfaces \cite{KN}.

This leads to the following description of the derived category: if
we denote the direct sum of the collection by $T = \oplus_{i=1}^r
T_i$ then $T$ is a tilting object, i.e.
\begin{equation}\label{tiltingequivalence}
\RHom(T,-) : D^b(Z) \iso D^b(\End_Z(T))
\end{equation}
is a triangulated equivalence. The fact that we use only $\End_Z(T)$
instead of $\RHom_Z(T,T)$ is because there are no higher Ext's (the
collection is strong), and the fact that this is an equivalence is
because the collection generates the whole derived category (it is
full). The astute reader will have noted
that this functor produces right modules not left modules, so to be consistent with the rest
of this paper we should (but won't) replace $\End_Z(T)$ with its opposite algebra.

We make one further assumption on $\set{T_i}$, that for any $i,j$
and $p>0$ we have
$$\Ext^k_Z(T_i, T_j\otimes \omega_Z^{-p}) = 0$$
for all $k>0$. The case $p=0$ is just the meaning of the word
\quotes{strong}. We shall call this a \textit{simple} collection (after \cite{bridgeland}),
it is a generalization of what Bondal
and Polishchuk call \quotes{geometric} \cite{BP}, though their definition is
in terms of mutations of the collection. However they show that a
geometric collection can only exist on a variety where
$$\mbox{rk } K(Z) = \mbox{dim} (Z) +1$$
in which case (as they also show) their definition is equivalent to
ours. Simple collections exist on all the del Pezzos and on the ruled surfaces
at least up to $\mathbb{F}_2$.

The algebra $A:=\End_Z(T)$ has a very simple form. Since $T$ is the
direct sum of $r$ objects it is clear that $A$ may be described
as the path algebra of a quiver with $r$ nodes plus some relations.
Also, by the axioms for a full strong exceptional collection, we may
order the $T_i$ so that the Hom's only go in one direction (say
increasing $i$) and so the quiver is directed. Thus we may give $A$
an $\mathbb{N}$-grading by declaring $\Hom_Z(T_i, T_j)$ to be of
degree $j-i$. The degree zero piece is just $\oplus_i \End_Z(T_i) =
\mathbb{C}^r$.
 Hence $A$ is of the correct form for Theorem
\ref{gradedresult} to apply, so it is given by a presentation
\begin{equation}\label{presentationforA}
m^\vee: \Ext_A^2(A_0,A_0)^\vee \to T \Ext_A^1(A_0,A_0)^\vee
\end{equation}
for $m$ a (graded) $A_\infty$-structure on $\Ext_A(A_0,A_0)$.

Now form the local Calabi-Yau $\omega_Z$. It is straight-forward
\cite{bridgeland} to show that the pull-ups of the $T_i$ to
$\omega_Z$ generate $D^b(\omega_Z)$. The projection formula gives
\begin{equation}\label{End_omegaT}
\Ext_\omega (\pi^*T, \pi^*T) = \bigoplus_{p\geq 0 } \Ext_Z(T,
T\otimes \omega_Z^{-p})
\end{equation}
so by our condition all the higher self-Ext's of $\pi^*T$ vanish.
Thus $\pi_* T$ is also a tilting object and
$$D^b(\omega) \cong D^b(\End_\omega(\pi^*T))$$
The question we asked in the introduction was: given a presentation
of the form (\ref{presentationforA}) for $A$, can we give a
presentation of
$$\tilde{A}:= \End_\omega(\pi_*T) ?$$
The answer, we claimed, is that $\tilde{A}$ is given by the cyclic
completion of the quiver corresponding to the presentation
(\ref{presentationforA}).

We now fill in the final details in the justification of this
answer. Using (\ref{End_omegaT}) we see we may also give $\tilde{A}$
an $\mathbb{N}$-grading, if we declare
$$\Hom_Z(T_i, T_j\otimes\omega_Z^{-p})$$
 to have degree $j-i + rp$. The degree-zero piece is still $\mathbb{C}^r$, so $\tilde{A}$
 also admits a presentation of the form
\begin{equation}\label{presentationfortildeA}m^\vee: \Ext_{\tilde{A}}^2(\tilde{A}_0,\tilde{A}_0)^\vee \to
 T
 \Ext_{\tilde{A}}^1(\tilde{A}_0,\tilde{A}_0)^\vee\end{equation}
Let
$$S=\bigoplus_{i=1}^r S_i \in D^b(Z)$$
denote the image of $A_0\in D^b(A)$ under the derived equivalence
(\ref{tiltingequivalence}). The $S_i$ form a dual exceptional
collection to the $T_i$. One checks \cite{bridgeland} that
$\iota_*S$ is the object in $D^b(\omega_Z)$ corresponding to
$\tilde{A}_0$. We claim (Lemma \ref{tiltingisdgequiv} below) that a
derived equivalence obtained by tilting is necessarily an
$A_\infty$-equivalence, so
$$\Ext_A(A_0,A_0) = \Ext_Z (S,S)$$
and
$$\Ext_{\tilde{A}}(\tilde{A}_0,\tilde{A}_0) = \Ext_\omega
(\iota_*S,\iota_* S)$$
 as $A_\infty$-algebras. Thus by Theorem
\ref{Extonomegaiscycliccompletion}
$\Ext_{\tilde{A}}(\tilde{A}_0,\tilde{A}_0)$ is the 3-dimensional
cyclic completion of $\Ext_A(A_0,A_0)$. This means that the
presentation (\ref{presentationfortildeA}) is really the map
$$\Ext_A^2(A_0,A_0)^\vee \oplus  \Ext_A^1(A_0,A_0) \to T\left( \Ext_A^1(A_0,A_0)^\vee
\oplus \Ext_A^2(A_0,A_0)\right)$$ given by dualising the cyclic completion of the
map $m$ in
(\ref{presentationforA}). This corresponds precisely to the process
we described of cyclicly completing the quiver.

This presentation
may be encoded in a superpotential using the construction from Theorem \ref{CY3impliesSP}.

There is one element missing in this story - really we should give a
criterion for an arbitrary presentation of $A$ to arise from an
$A_\infty$ structure in the manner of (\ref{presentationforA}). It
seems plausible that the proof in \cite{LPWZ} might yield such a
criterion. As it is we shall just assume that any reasonable
presentation does arise in this way.

\begin{lem}\label{tiltingisdgequiv} Suppose that the derived category $D^b(\mathcal{C})$ of some abelian category
$\mathcal{C}$ has a tilting object $T = \oplus T_i$, and let
$A=\End_Z(T)$. Then the derived equivalence
$$\Psi = \RHom_\mathcal{C}(T,-) : D^b(\mathcal{C}) \iso D^b(A)$$
is in fact a dg (or $A_\infty$) equivalence, i.e. for any $E\in
\mathcal{C}$ we have
$$\RHom_\mathcal{C}(E,E) \simeq \RHom_A(\Psi E, \Psi E)$$
as dgas.
\end{lem}
\begin{proof}
Sending $E$ through $\Psi$ and back produces a
resolution $\mathcal{E} \simeq E$ in terms of the $T_i$. There are no higher
$\Ext$s between the $T_i$, so
$$\RHom_\mathcal{C}(E,E)=\Hom_\mathcal{C}(\mathcal{E},\mathcal{E})$$
For the same reason $\Psi E = \Hom_\mathcal{C}(T,\mathcal{E})$, which is a complex $\mathcal{F}$ of the projective modules $\set{A_i:=\Psi T_i}$ that is isomorphic to $\mathcal{E}$.
Thus
$$\RHom_A(\Psi E, \Psi E) = \Hom_A(\mathcal{F},\mathcal{F}) = \Hom_\mathcal{C}(\mathcal{E},\mathcal{E})$$
\end{proof}


\providecommand{\bysame}{\leavevmode\hbox
to3em{\hrulefill}\thinspace}
\providecommand{\MR}{\relax\ifhmode\unskip\space\fi MR }
\providecommand{\MRhref}[2]{%
  \href{http://www.ams.org/mathscinet-getitem?mr=#1}{#2}
} \providecommand{\href}[2]{#2}

\end{document}